\documentclass[a4,amsmath, amsthm, amssymb]{amsart}

\usepackage{pdfpages}
\usepackage[left=20mm,top=0.6in,bottom=12mm]{geometry}

\numberwithin{equation}{section}

\newtheorem{Lemma}{LEMMA}[section]
\newtheorem{Theorem}[Lemma]{Theorem}
\newtheorem{Proposition}[Lemma]{Proposition}
\newtheorem{Corollary}[Lemma]{Corollary}

\newtheorem{remark}[Lemma]{Remark}
\newtheorem{definition}[Lemma]{Definition}
\newtheorem{example}[Lemma]{Example}

\newtheorem{Fact}[Lemma]{Fact}
\newtheorem{assumption}[Lemma]{Assumption}
\def\bt{\begin{Theorem}}
\def\et{\end{Theorem}}
\def\bl{\begin{Lemma}}
\def\el{\end{Lemma}}
\def\bp{\begin{Proposition}}
\def\ep{\end{Proposition}}
\def\bcor{\begin{Corollary}}
\def\ecor{\end{Corollary}}
\def\bpf{\begin{proof}}
\def\epf{\end{proof}}

\def\brem{\begin{remark}\rm }
\def\erem{\end{remark}}

\def\bedef{\begin{definition}\rm }
\def\endef{\end{definition}}

\def\beg{\begin{example}}
\def\eeg{\end{example}}

\def\bef{\begin{Fact}}
\def\eef{\end{Fact}}

\def\bea{\begin{assumption}}
\def\ena{\end{assumption}}
\def\bc{\begin{center}}
\def\ec{\end{center}}
\def\noi{\noindent}

\def\vsq{\vskip .25cm}
\def\beq{\begin{equation}}
\def\eeq{\end{equation}}
\def\beqarray{\begin{eqnarray*}}
\def\eeqarray{\end{eqnarray*}}
\def\<{\leftangle}
\def\>{\rightangle}
\def\({\left(}
\def\){\right)}
\def\f{\varphi}

\def\<{\langle}
\def\>{\rangle}
\def\q{\quad}
\def\r{\rho}
\def\a{\alpha}
\def\b{\beta}
\def\g{\gamma}

\def\d{\delta}
\def\h{\hbox}

\def\t{\tau}
\def\r{\rho}
\def\l{\lambda}
\def\e{\varepsilon}

\def\O{\Omega}

\def\B{$\blacksquare$}
\def\w.r.t.{with respect to}
\def\R{{\mathbb{R}}}
\def\N{{\mathbb{N}}}

\def\C{{\mathbb{C}}}

\def\bq{\begin{quote}}
\def\eq{\end{quote}}

\def\bit{\begin{itemize}}
\def\eit{\end{itemize}}
\def\i{\item}
\def\ben{\begin{enumerate}}
\def\een{\end{enumerate}}

\def\ds{\displaystyle}

\def\X{{\mathcal X}}

\def\A{{\mathbf A}}
\def\B{{\mathbf B}}
\def\C{{\mathbf{C}}}
\def\D{{\mathbf D}}
\def\T{{\mathcal T}}
\def\P{{\mathbf P}}

\begin{document}
\title[Identification of matrix coefficients]{ Identification of matrix diffusion coefficients\\ in a parabolic PDE }
\author{Subhankar Mondal and M. Thamban Nair}
\address{Department of Mathematics, IIT Madras, Chennai 600036, INDIA}
\email{s.subhankar80@gmail.com; mtnair@iitm.ac.in}
\maketitle

%\begin{abstract}We consider an inverse problem of identifying the diffusion coefficient in  matrix form in a parabolic PDE. In 2006, Cao and Pereverzev, used a \textit{natural linearisation} method for identifying a scalar valued diffusion coefficient in a parabolic PDE. In this paper, we make use of that idea for identifying a matrix valued coefficient. Under some assumption, we show the uniqueness of the inverse problem. Using the notion of a weak solution for a parabolic PDE, we transform our non-linear inverse problem into a problem of solving an operator equation where the operator involved is linear, and which turns out to be ill-posed. An explicit representation of adjoint of the linear operator involved is obtained. The method of Tikhonov regularization is employed for obtaining stable approximations and for finite dimensional analysis the Galerkin method is used. For the analysis in the finite dimensional setting, we have defined orthogonal projections on the space of matrices with entries from $L^2(\O).$ For choosing the regularizing parameter, we used the adaptive technique, so that we have an order optimal rate of convergence. Finally, for the relaxed noisy data,  we described a procedure for obtaining a smoothed version so as to obtain the error estimates. 
%\end{abstract}

\begin{abstract}We consider an inverse problem of identifying the diffusion coefficient in  matrix form in a parabolic PDE. In 2006, Cao and Pereverzev, used a \textit{natural linearisation} method for identifying a scalar valued diffusion coefficient in a parabolic PDE. In this paper, we make use of that idea for identifying a matrix valued coefficient, namely, using the notion of a weak solution for a parabolic PDE, we transform our non-linear inverse problem into a problem of solving an ill-posed  operator equation where the operator depending on the data  is linear. For the purpose of obtaining stable  approximate solutions, Tikhonov regularization is employed, and error estimates under noisy data are derived.  We have also showed the uniqueness of the solution of the inverse problem under  some assumptions on the data and obtained  explicit representation of adjoint of the linear operator involved. For the obtaining error estimates in the  finite dimensional setting,  Galerkin method is used, by defining  orthogonal projections on the space of matrices with entries from $L^2(\O)$, by making use of standard orthogonal projections  on $L^2(\O)$. For choosing the regularizing parameter, we used the adaptive technique, so that we have an order optimal rate of convergence. Finally, for the relaxed noisy data,  we described a procedure for obtaining a smoothed version so as to obtain the error estimates. 
\end{abstract}



\section{Introduction}\label{intro}
Let $\O$ be a bounded domain in $\R^d$ with Lipschitz boundary $\partial \O.$ For a fixed $\tau >0$, we denote $\O\times (0, \tau)$ by $Q$ and its boundary $\partial\O\times (0,\tau)$ by $\partial Q,$ also we denote the interval $[0,\tau]$ by $I_\tau.$ We consider an inverse problem associated with the PDE
\beq\label{pde1}
u_t-\nabla\cdot (\A(x)\nabla u)=f(x,t) \q \text{in}\ Q
\eeq
along with the following boundary and initial conditions:
\beq\label{bdry}
u(x, t)= 0\q \text{on}\q \partial Q,
\eeq
\beq\label{initial}
u(x, 0)=h(x)\ \text{in}\ \O,
\eeq
where 
% {\cred $\A\in  \big(H^1(\O)\big)^{d\times d},$ that is, 
{$\A$ is a $d\times d$ matrix with entries from} $H^1(\O)$, $h\in L^2(\O)$ and $f\in L^2(0, \tau; L^2(\O))$.

{In the above,  and throughout the paper,
if $\mathcal{H}$ is a Hilbert space, then  by $L^2(0, \tau; \mathcal{H})$ we mean the space  $L^2(I_\tau, {\mathcal H})$, that is, 
the space of all measurable functions  $\psi: [0, \tau]\to \mathcal{H}$, defined a.e. on $[0, \t]$,  such 
that 
$$\int_0^\t \|\psi(t)\|_{\mathcal{H}}^2 dt <\infty.$$
Also, if $\mathcal{H}$ is a Hilbert space of functions or equivalence classes of functions defined a.e. on $\O$ and if $\psi: \O\times [0, \tau]\to \R$ is an a.e. defined measurable function, then we  write $\psi\in L^2(0, \tau; \mathcal{H})$ if and only if $\tilde\psi: [0, \tau]\to {\mathcal H}$ defined by
$$\tilde\psi(t)(x) = \psi(x, t)\q\h{for almost all}\q (x, t)\in \O\times [0, \tau],$$
belongs to $L^2(I_\tau, {\mathcal H})$.}

\noindent
%{\cred Throughout the paper, $\mathcal{H}$ can be $H^1_0(\O), L^2(\O)$ or $H^{-1}(\O),$ the dual of $H^1_0(\O).$}

Using the above notation, we may observe that
$$f\in L^2(0, \tau;  L^2(\O)) \iff f\in L^2(\O\times I_\tau).$$

\noi
{{\bf Convention:} 
%For $\psi \in L^2(0,\tau;\mathcal{H}(\O))$, we use the standard notation, $\nabla \psi$ 
{For $\psi \in L^2(0,\tau; H^1(\O))$, we use the notation, $\nabla \psi$ }for the gradient of $\psi$ with respect to the space variable. %{\color{blue}whenever it exists}.}

\noindent{
{\bf Notation:} For the Sobolev space $H^1_0(\O)$, its dual space will be denoted by standard notation $H^{-1}(\O)$ and for simplicity of the notations, their duality action will be denoted by $\<\cdot,\cdot\>$, that is for $\psi\in H^{-1}(\O)$ and $\f\in H^1_0(\O),$ $\<\psi,\f\>=\psi(\f).$}

\noindent The usual forward problem associated with \eqref{pde1}-\eqref{initial} is to find $u$ in some suitable space so that it satisfies \eqref{pde1}-\eqref{initial}. Most often one may be looking for a $u$ satisfying a weak formulation of \eqref{pde1}-\eqref{initial}. 
%In the next section we shall specify such a weak formulation.
The weak solution of \eqref{pde1}-\eqref{initial} is as per the following definition (cf. \cite{evans}).

{
\bedef\label{def-weakdef}
Let $\A,\,  f$ and $h$ be as in \eqref{pde1}-\eqref{initial}.  Then  $u\in L^2(0, \tau;  H^1_0(\O))$ with $u_t \in L^2(0, \tau;  H^{-1}(\O))$ is said to be a {\bf weak solution} of the  parabolic system \eqref{pde1}-\eqref{initial} if $u$ satisfies 
\beq\label{weakform-1}
\langle u_t(\cdot ,t), \varphi\rangle+ \int_\O \A(x)\nabla u(\cdot, t)\cdot \nabla \varphi = \int_\O f(\cdot,t)\varphi
\eeq
for all $\varphi \in H_0^1(\O)$ and for almost all $t\in [0,\tau]$,
along with $u(x, 0)=h(x)$ for  $x\in \O.$ 
If $f$ belongs to $L^2(0, \tau;  H^{-1}(\O))$  instead of $L^2(0, \tau;  L^2(\O))$, then the R.H.S of  \eqref{weakform-1} may be replaced by the duality pairing $\langle f(.,t),\varphi\rangle$.
\endef
}

{
A natural question is: Under what conditions on $\A,\, f, \, h$ we can ensure the existence of a weak solution for  the  parabolic system \eqref{pde1}-\eqref{initial}?  

The  theorem below (Thorem \ref{existence-evans}) specifies certain conditions under which the above  question is answered affirmatively (see Theorem 5, Theorem 6 in Chapter 7 of \cite{evans}).  Before stating the theorem,  let us recall the following standard definitions: 

 \bedef\label{def-symmetric-elliptic}  Let  $\B$ be a  $d\times d$ matrix with entries from $L^2(\O)$. 
 \ben
 \i[(i)] $\B$ is  said to be {\bf symmetric} a.e. on $\O$ if the matrix $\B(x)$ is symmetric for almost all $x\in \O$.

\i[(ii)] $\B$ is said to satisfy the {\bf uniform ellipticity condition}, if there exists a $q_0>0$ such that
%\beq\label{elliptic}
$$\B(x)\xi\cdot\xi=\sum_{1\leq i,j\leq d}B_{ij}(x)\xi_i\xi_j\geq q_0 |\xi |^2 \h{ for almost all }  x\in \O$$
%\eeq
and for all $\xi=(\xi_1,...,\xi_d)\in \R^d$, where $|\xi|= (\xi_1^2+\xi_2^2+...+\xi_d^2)^{1/2}$.
\een
\endef
}

{\bt\label{existence-evans}\cite[Theorem 6, Ch. 7]{evans}
Let $\A$ be a $d\times d$ symmetric matrix with entries from $H^1(\O)$ such that it satisfies the uniform ellipticity condition as in Definition \ref{def-symmetric-elliptic}.  Then \eqref{pde1}-\eqref{initial} has a unique {\it weak solution} $u\in L^2(0,\tau;H^{2m+2}(\O))$ provided  $h$ and $f$ are such that 
$h\in H^{2m+1}(\O)$ and 
$$\frac{\partial^kf}{\partial t^k}\in L^2(0,\tau;H^{2m-2k}(\O)) \q\h{for}\q  k=0,\ldots, m,\q\h{for some}\q  m\in \N\cup\{0\},$$
 and $h_0,\, h_1, \ldots, h_m$, defined iteratively as follows,  belong to $H^1_0(\O)$: 
$$h_0:=h,\q  h_1:=f(\cdot,0)-\nabla\cdot \A \nabla\,h_0,\q \cdots,\q  h_m:=\frac{\partial^{m-1}f(\cdot,0)}{\partial t^{m-1}}-\nabla\cdot\A\nabla\,h_{m-1}. $$
\et }

\noi\,
In this paper, we are  interested in finding regularized approximations for the   following inverse problem {({\bf IP})} associated with \eqref{pde1}-\eqref{initial}:% , denoted by (IP):
\bq
{\bf (IP):}\, To  identify the matrix diffusion coefficient $\A$ from the knowledge of $u$, a weak solution of \eqref{pde1}-\eqref{initial},  which may be known only approximately with some noise.
\eq

Clearly, the { above inverse problem ({\bf IP})  is non linear}. 
{We shall represent  this non linear inverse problem as  an  equation involving a linear operator and carry out the regularization procedure as in the case of a linear operator equation. We assume that the datum $u$ is such that this inverse problem ({\bf IP}) has a solution.  More precisely, in view of the Definition \ref{def-weakdef}  and Theorem \ref{existence-evans}, we make the following assumption. }

{
\bea\label{existence}  There exists a  $d\times d$  matrix  $\A$ with entries from {$H^1(\O)$}  such that  the  parabolic problem  \eqref{pde1}-\eqref{initial} has a  unique weak solution $u\in L^2(0, \tau; H^1_0(\O))$ with $u_t\in L^2(0,\tau;H^{-1}(\O))$ and $0\neq\nabla u\in L^2(0,\tau;L^\infty(\O,\R^d))$.  
\ena
}

{
\brem\, It is to be observed that the additional condition assumed for $u$ in Assumption \ref{existence}, namely,   $\nabla u\in L^2(0,\tau;L^\infty(\O,\R^d))$ is satisfied if $m>\frac{d}{2}$, and if $f$ and $h$ satisfy the hypothesis of Theorem \ref{existence-evans}, for in this case,  $\nabla u(\cdot,t)\in H^m(\O)$ and $H^m(\O)$ is continuously embedded in $L^\infty(\O)$ (cf. \cite[Theorem 2.5.4]{kesh}). 
\erem
}

Parameter identification problem has a vast literature and every problem has its own importance in practical application. One can refer
to \cite{engl, cao, hanke, tommi, yuan} for more literature on parameter identification problem in parabolic setting.
In \cite{cao}, the authors have considered a method called {\it natural linearisation} for identifying a scalar valued diffusion coefficient from final time observation. There, they have suitably reduced the problem involved with parabolic PDE into an elliptic PDE, with some assumptions, and then applied the natural linearisation technique for
identifying the parameter.

 In this paper, we use similar idea for converting the ({\bf IP}) into a problem of  solving an ill-posed operator equation, where the operator involved is linear, but without converting into an elliptic PDE. We will do the analysis with parabolic PDE setting. Also, we would like to mention that in \cite{cao},
the authors have worked with Neuman boundary condition, but in this work, we are considering homogeneous Dirichlet boundary condition.
Further, we are considering the problem in a more general setting, namely, the identification is that of a $d\times d$ matrix coefficient with functions from  {$H^1(\O)$} as entries, instead of scalar valued function.  We use the Tikhonov regularization method to obtain stable approximations from the space of all $d\times d$ matrix with entries from $L^2(\O).$
As a result, while doing analysis in the finite dimensional setting, we found it necessary to define an orthogonal projection on the space of all $d\times d$ matrices with entries from $L^2(\O)$. Also, we obtain an explicit representation of the adjoint of the linear operator involved (see Theorem \ref{adj-rep-2}).

For obtaining error estimates for the regularized approximations, first we assume certain regularity conditions on the noisy data. In order to accommodate the case with relaxed noisy data, we considered a smoothing procedure using  the {\it Clement operator}, which resulted in a reduced accuracy in the error.

{
We would like to point out that most of the parameter identification problems for parabolic PDE's are on identifying scalar valued coefficient functions, without paying much attention to the problem of identifying matrix valued coefficient functions. To the best of our knowledge, \cite{yuan} is the first paper that deals with uniqueness and stability for the inverse problem of identifying a matrix diffusion coefficient in a parabolic PDE from some boundary observations and an intermediate time observation, and not much has been available on  identification of matrix coefficients. Thus, this paper intents to make some contribution in that direction. Using some idea employed in \cite{yuan}, we could also establish a uniqueness result for our inverse problem under some additional condition on $u$ (see Theorem \ref{uniq-th-ip-2}).}

The paper is organized as follows. In Section 2, we introduce the {notations that will be used  throughout the paper and present  some preliminary results that will be required for further analysis, and also establish the uniqueness result for our inverse problem}. In Section 3, we reformulate the inverse problem as {an operator equation involving a linear operator} and prove some results related to boundedness, compactness and rank of the linear operator of our interest. Section 4 deals with regularization and error estimates. We use {the standard theory of  Tikhonov regularization} for regularization purpose. In Section 5 we considered the finite dimensional realizations of the regularized solutions and derived corresponding error estimates. Section 6 deals with the strategy of choosing regularizing parameter, based on the adaptive technique. 
Section 7 is devoted to the procedure of obtaining smoothed version of the observed data, and subsequently, in Section 8 we have given the error estimates for the smoothed version of the noisy data.
%\newpage

\section{Preliminaries and a uniqueness result}\label{prelim}
We shall consider the real vector space
\beq\label{space-new} \X=\{v\in L^2(0, \tau;  H^1_0(\O)):v_t\in L^2(0, \tau;  H^{-1}(\O))\}.\eeq
 It can be seen that
$$\lVert{v}\rVert_{\X}=\lVert{v}\rVert_{L^2(0, \tau;  H^1_0(\O))}+\lVert{v_t}\rVert_{L^2(0, \tau;  H^{-1}(\O))},\q  v\in \X,$$
defines a norm on $\X$ which makes it a Banach space.

{Throughout  the paper, we shall use the notation ${\mathcal H}^{d\times d}$ to denote the space of all $d\times d$ matrices with entries from a Hilbert space ${\mathcal H}$.  }
Also,  we consider the inner product { 
$$\Big\< \B,\tilde{\B}\Big\>=\sum_{i,j=1}^d\langle b_{ij}, \tilde{b}_{ij}\rangle_{{L^2(\O)}},$$
and the corresponding norm
$$\|\B\|= \Big(\sum_{i,j=1}^d \|b_{ij}\|^2_{{L^2(\O)}}\Big)^{1/2}$$
for  $\B=(b_{ij})$ and $\tilde \B = (\tilde b_{ij})$ in  ${(L^2(\O))}^{d\times d}$}.
%, where $F$ in $\|\cdot\|_F$ stands for {\it Frobenius}. 

{
\noindent
{\bf Notations:} For simplicity of the presentation we shall use the notations $H^1$, $W^{1,\infty}$, $L^2$, $L^\infty$ for denoting the spaces $H^1(\O)$, $W^{1,\infty}(\O)$, $L^2(\O)$ and $L^\infty(\O),$ respectively. 
}
\vsq

Let $1\leq p\leq \infty.$  For $\vec{v}: \O\to \R^d$, if  $\vec{v}(x)=(v_1(x),...,v_d(x)),\, x\in \O$, then we see that
$$\vec{v}\in L^p(\O,\R^d) \iff v_i\in L^p(\O)\q\forall i\in \{1, \ldots, d\}$$
and
$$\ds\lVert{\vec{v}}\rVert^p_{L^p}=  \left\{\begin{array}{ll}
                                      \sum_{i=1}^d\lVert{v_i}\rVert^p_{L^p}&\h{for}\q 1\leq p<\infty,\\
                                       \max\{\lVert{v_i}\rVert_{L^\infty}: 1\leq i\leq d\} &\h{for}\q p=\infty. \end{array}\right.$$
                                       
Let $\A_0\in (L^2(\O))^{d\times d}$ be symmetric and satisfies the uniform ellipticity condition (see Definition \ref{def-symmetric-elliptic}). For $f\in L^2(0, \tau;  L^2(\O))$ and $h\in L^2(\O)$,  consider the PDE
\beq\label{pde2}
v_t-\nabla\cdot (\A_0(x)\nabla v)=f(x,t) \q \text{in}\ Q
\eeq
along with the boundary and initial conditions
\beq\label{bdry2}
v(x, t)= 0\q \text{on}\ \partial Q,
\eeq
\beq\label{initial2}
v(x, 0)=h(x)\q \text{in}\ \O.
\eeq
The equation (\ref{bdry2}) is to be understood in the sense of {\it trace}.
Then we have the following result on  existence and uniqueness.

\bt\label{dirichletestimate} \cite[Theorem 2 - 4, Ch. 7]{evans}
{ Let $\A_0, f$ and $h$ be as in \eqref{pde2}-\eqref{initial2}. 
Then \eqref{pde2}-\eqref{initial2} has a  unique $v\in L^2(0, \tau;  H^1_0(\O))$  with $v_t\in L^2(0, \tau;  H^{-1}(\O))$ satisfying
\beq\label{weakform}
\langle v_t(\cdot ,t), \varphi\rangle + \int_\O \A_0(x)\nabla v(\cdot, t)\cdot \nabla \varphi = \int_\O f(\cdot,t)\varphi
\eeq
{and $v(x, 0)=h(x)$ for  $x\in \O,$}
for all $\varphi \in H_0^1(\O)$ and for almost all $t\in [0,\tau]$.}
 Further, there exist a constant $C_0>0$ depending on $\O, \tau$ and $\A_0$, such that
\beq\label{dirichletestimatev}
\lVert{v}\rVert_{L^2(0, \tau;  H^1_0(\O))}+\lVert{v_t}\rVert_{L^2(0, \tau;  H^{-1}(\O))}\leq C_0\Big(\lVert{f}\rVert_{L^2(0, \tau;  H^{-1}(\O))}+\lVert{h}\rVert_{L^2(\O)}\Big).
\eeq
\et

In this  paper we are dealing with spaces of Hilbert space valued functions, for instance $L^2(0, \tau;  H^1_0(\O)),$ $ L^2(0, \tau;  L^2(\O)).$ So we would like to have some compact embeddings for these type of function spaces, which will be useful { in proving some convergence results. } With regard to this,  we have the following result, known as {\bf Aubin-Lions} lemma which gives a compact embedding between certain Banach space valued function spaces. For its proof and more details about Aubin-Lions lemma one may refer to any of  \cite{aubin, serano, simons}.

Let $X_0$ and $X_1$ be Banach spaces with $X_0\subset X_1$ and
\beq\label{w}\mathcal{W}=\{u\in L^2(0, \tau;  X_0): u_t\in L^2(0, \tau;  X_1)\}.\eeq Then it can be seen that $\mathcal{W}$ is a Banach space with respect to the norm
$$\lVert{u}\rVert_{\mathcal{W}}=\lVert{u}\rVert_{L^2(0, \tau;  X_0)}+\lVert{u_t}\rVert_{L^2(0, \tau;  X_1)}, \q u\in {\mathcal W}.$$

\bl\text{{\bf (Aubin-Lions)}}\label{aubin-lion}\cite[Theorem 1.3]{serano} 
Let $X_0,X,X_1$ be Banach spaces with $X_0\subset X\subset X_1.$ Let $\mathcal{W}$ be the Banach space as defined in \eqref{w}.
If $X_0$ is compactly embedded in $X$ and $X$ is continuously embedded in $X_1$, then  $\mathcal{W}$ is compactly embedded in $L^2(0, \tau;  X)$.
\el

\brem\label{compactembedding}
Let $\X$ be the space as in \eqref{space-new}. Under our assumption on $\O$, $H^1_0(\O)$ is compactly embedded in $L^2(\O)$ and $L^2(\O)$ is continuously embedded in $H^{-1}(\O)$ (see \cite{kesh}). Therefore, by Lemma \ref{aubin-lion}, $\X$ is compactly embedded in $L^2(0, \tau;  L^2(\O)).$
\erem

We know that if $\psi\in L^2(\O)$, then $\dfrac{\partial \psi}{\partial x_i}$,  in the sense of distribution, belongs to $H^{-1}(\O)$ for all $1\leq i \leq d$ (see \cite{kesh}). Next, we have a simple but useful result.

\bl\label{h-1bound}
If $\psi\in L^2(\O)$, then $\left\lVert{\dfrac{\partial \psi}{\partial x_i}}\right\rVert_{H^{-1}}\leq \lVert{\psi}\rVert_{L^2}$  for each $1\leq i\leq d$.
\el

\bpf
Let $\varphi\in C_c^\infty(\O),$ the set of all real valued infinitely differentiable functions on $\O$ with compact support. Then, for $1\leq i\leq d$, we have $$\big(\dfrac{\partial \psi}{\partial x_i}\big)(\varphi)=-\int_\O\psi \dfrac{\partial \varphi}{\partial x_i}\q \text{for all}\ \varphi\in C^\infty_c(\O). $$
 Thus, for all $\varphi\in C_c^\infty(\O)$
$$
 \left |\dfrac{\partial \psi}{\partial x_i}(\varphi)\right | \leq  \int_\O\ |\psi \dfrac{\partial \varphi}{\partial x_i} | \leq   \lVert{\psi}\rVert_{L^2}\left\lVert{\dfrac{\partial \varphi}{\partial x_i}}\right\rVert_{L^2}
\leq  \lVert{\psi}\rVert_{L^2}\lVert{\varphi}\rVert_{H^1(\O)}.$$
Since $C_c^\infty(\O)$ is dense in $H^1_0(\O)$, we have $\|{\frac{\partial \psi}{\partial x_i}}\|_{H^{-1}}\leq \lVert{\psi}\rVert_{L^2}$ for $1\leq i\leq d$.
\epf

\brem\label{vectorh-1bound}
Let $\psi\in L^2(\O).$ From Lemma \ref{h-1bound}, we also have $$\lVert{\nabla \psi}\rVert^2_{H^{-1}}= \sum_{i=1}^d\left\lVert{\dfrac{\partial \psi}{\partial x_i}}\right\rVert^2_{H^{-1}}\leq d \lVert{\psi}\rVert^2_{L^2}.$$ That is, $\lVert{\nabla \psi}\rVert_{H^{-1}}\leq \sqrt{d}\lVert{\psi}\rVert_{L^2}.$
\erem

\subsection{Uniqueness of the inverse problem}  {We establish the uniqueness result for our inverse problem {\bf (IP)} under appropriate assumptions. {The arguments used here are  motivated by those in  \cite{yuan}.}}

Let $\A\in (H^1(\O))^{d\times d}$ be a solution of the inverse problem {\bf (IP)} corresponding to the data $u$. Then, we have 
\beq\label{uniq-ip-2}
u_t-\nabla\cdot\A\nabla u=f\q\text{in}\,\,\O\times (0,\tau).
\eeq
Suppose $\B\in (H^1(\O))^{d\times d}$ is another solution to {\bf (IP)}. Then, from \eqref{uniq-ip-2}, we have 
\beq\label{red-uniq-ip-2}
-\nabla\cdot (\A-\B)\nabla u=0\q\text{in}\,\,\O\times (0,\tau).
\eeq
 {Let $t_i\in (0, \t)$ for $ i=1, \ldots, d^2(d+1)$.  We denote $u(\cdot,t_l)$  and $\frac{\partial u(t_l)}{\partial x_i}$  by $u(t_l)$ and  $\partial_iu(t_l)$, respectively, and for  $1\leq k\leq d^2$ and $1\leq i,j\leq d,$ we let}
$$D^k_{ij}=\det\begin{pmatrix}
\partial_1u(t_{(k-1)(d+1)+1})&\cdots& \partial_d u(t_{(k-1)(d+1)+1})&\partial_i\partial_j u(t_{(k-1)(d+1)+1})\\
\cdot&\cdots&\cdot&\cdot\\
\cdot&\cdots&\cdot&\cdot\\
\partial_1u(t_{(k-1)(d+1)+d+1})&\cdots&\partial_du(t_{(k-1)(d+1)+d+1})&\partial_i\partial_ju(t_{(k-1)(d+1)+d+1})
\end{pmatrix}
$$
and \beq\label{det-2}
D=\det\begin{pmatrix}
D^1_{11}&\cdots &D^1_{1d}&D^1_{21}&\cdots&D^1_{2d}&\cdots &D^1_{d1}&\cdots&D^1_{dd}\\
\cdot&\cdots&\cdot&\cdot&\cdots&\cdot&\cdots&\cdot&\cdots&\cdot\\
\cdot&\cdots&\cdot&\cdot&\cdots&\cdot&\cdots&\cdot&\cdots&\cdot\\
D^{d^2}_{11}&\cdots&D^{d^2}_{1d}&D^{d^2}_{21}&\cdots&D^{d^2}_{2d}&\cdots&D^{d^2}_{d1}&\cdots&D^{d^2}_{dd}
\end{pmatrix}.
\eeq

\bt\label{uniq-th-ip-2}
Let $\Big\{t_i\,:\,1\leq i\leq d^2(d+1)\Big\}\subset (0,\tau)$ and $\A,\B\in (H^1(\O))^{d\times d}$ {be  solutions} of the inverse problem \textbf{(IP)} corresponding to the data $u\in \X$. Let $D$ be as defined in \eqref{det-2}. If $D\neq 0$ a.e. in $\O$, then $\A=\B$ a.e. in $\O.$
\et
\bpf
From the assumptions on $\A$ and $\B$, it follows that they satisfy \eqref{red-uniq-ip-2}.
Let $f_{ij}$ be the $ij$-th entry of the matrix $(\A-\B)$. Thus, for $1\leq l\leq d^2(d+1),$ from \eqref{red-uniq-ip-2}, we have
\beq\label{system-eq-2}
\Big(\sum_{j=1}^d\partial_1f_{1j}\partial_ju(t_l)+\cdots+\sum_{j=1}^d\partial_df_{dj}\partial_ju(t_l)\Big)=-\Big(\sum_{i,j=1}^df_{ij}\partial_i\partial_ju(t_l)\Big).
\eeq
Therefore, for $1\leq k\leq d^2,$ from \eqref{system-eq-2}, we have
\beqarray\begin{pmatrix}
\partial_1u(t_{(k-1)(d+1)+1})&\cdots&\partial_du(t_{(k-1)(d+1)+1})\\
\cdot&\cdots&\cdot\\
\cdot&\cdots&\cdot\\
\partial_1u(t_{(k-1)(d+1)+d+1})&\cdots&\partial_du(t_{(k-1)(d+1)+d+1})
\end{pmatrix}&&
\begin{pmatrix}
\sum_{j=1}^d\partial_jf_{j1}\\
\cdot\\
\cdot\\
\sum_{j=1}^d \partial_jf_{jd}
\end{pmatrix}\\
&&=\begin{pmatrix}
-\sum_{i,j=1}^df_{ij}\partial_i\partial_ju(t_{(k-1)(d+1)+1})\\
\cdot\\
\cdot\\
-\sum_{i,j=1}^df_{ij}\partial_i\partial_ju(t_{(k-1)(d+1)+d+1})
\end{pmatrix}
\eeqarray
For each $1\leq k\leq d^2,$ the above system of $d+1$ equations  in $d$ variables has a solution, namely $$\Big(\sum_{i,j=1}^d\partial_jf_{j1},\cdots,\sum_{i,j=1}^d\partial_jf_{jd}\Big)^T,$$
therefore, we must have 
$$\det\begin{pmatrix}
\partial_1u(t_{(k-1)(d+1)+1})&\cdots&\partial_du(t_{(k-1)(d+1)+1})&-\sum_{i,j=1}^df_{ij}\partial_i\partial_ju(t_{(k-1)(d+1)+1})\\
\cdot&\cdots&\cdot&\cdot\\
\cdot&\cdots&\cdot&\cdot\\
\partial_1u(t_{(k-1)(d+1)+d+1})&\cdots&\partial_du(t_{(k-1)(d+1)+d+1})&-\sum_{i,j=1}^df_{ij}\partial_i\partial_ju(t_{(k-1)(d+1)+d+1})
\end{pmatrix}=0.$$
Now, using the properties of determinant, we obtain the system of $d^2$ homogeneous equations with $d^2$ variables, namely,
\beq\label{hom-system-eq-2}
\sum_{i,j=1}^dD^k_{ij}f_{ij}=0,\q\text{for all}\,\, 1\leq k\leq d^2.
\eeq
By the hypothesis $D\neq 0$ a.e. in $\O,$ that is, determinant of the coefficient matrix of the system \eqref{hom-system-eq-2} is non-zero. Thus, $f_{ij}=0$ a.e. in $\O$ and hence the proof is complete.
\epf

In view of Theorem \ref{uniq-th-ip-2},    to guarantee the unique solvability of our proposed inverse problem, we also make the  one more assumption: 
 
{
 \bea\label{det-non-zero}
With   $u$ as in Assumption \ref{existence}, there exists $t_i\in (0, \t)$ with $1\leq i\leq d^2(d+1)$ such that $D\neq 0$ a.e. in $\O$, $D$ is as in (\ref{det-2}).
\ena
}

\section{Operator theoretic formulation}
In \cite{cao}, the authors have used the technique of natural linearisation for identification of a scalar valued diffusion coefficient from final time observation. In this section we show that  similar idea, {with some modifications in the arguments,  can be used for representing  our non linear ill-posed problem in an equivalent form that would facilitate} {us to make use of the theory of regularization for linear operator, as shown in the next section. Essentially, the method of natural linearisation, allows us to obtain solution of the non linear inverse problem in terms of solution of an operator equation, where the operator involved is a linear operator. This fact can be observed as we reach towards the end of this section.}

{
According to our inverse problem, $u$ satisfies \eqref{pde1}-\eqref{initial} in the weak sense. Thus, we have
\beq
\label{inv-weak} 
\langle u_t(\cdot ,t), \varphi\rangle + \int_\O \A(x)\nabla u(\cdot,t)\cdot \nabla \varphi = \int_\O f(\cdot,t)\varphi
\eeq
for all $\ \varphi \in H_0^1(\O)$ {and}  a.a. $t\in [0,\tau]$
along with
$u(x,0) = h(x)$  in  $\O.$ We shall consider an equivalent form of the above equation by making use of another matrix  $\A_0$, considered as in Theorem \ref{dirichletestimate}.

So, let  $\A_0$ be a  symmetric $d\times d$-matrix with entries from $L^2(\O)$  and  satisfies the uniform ellipticity condition. 
Let $f,\ h$ be as considered in the inverse problem {\bf (IP)}. By Theorem \ref{dirichletestimate}, there exists a unique weak solution   $v_0\in \X$ for \eqref{pde2}-\eqref{initial2}. That is, we have
\beq\label{v0}
\langle (v_0)_t(\cdot ,t), \varphi\rangle + \int_\O \A_0 \nabla v_0(\cdot, t)  \cdot \nabla \varphi  = \int_\O f(\cdot,t)\varphi \eeq
for all $\varphi \in H_0^1(\O)$  and  for a.a. $t\in [0,\tau]$
along with
$v_0(x,0)= h(x)$  in $\O.$   Thus, from (\ref{inv-weak}) and (\ref{v0}),  we have
}
\beq\label{v0-u}
\langle(v_0-u)_t(\cdot,t),\varphi\rangle+\int_\O \A_0\nabla (v_0-u)(\cdot,t)\cdot \nabla \varphi=\int_\O (\A-\A_0)\nabla u(\cdot, t)\cdot\nabla \varphi
\eeq for all $\varphi \in H_0^1(\O)$ along with $(v_0-u)(\cdot,0)=0$ in $\O$ and $ (v_0-u)(\cdot,t)=0$ on $\partial\O$ for a.a. $t\in [0,\tau].$
%We now state and prove a result that will ensure the existence of a weak solution of a PDE, that we will consider later.

{The following lemma will be used to show the existence of weak solution of a PDE that arises out in the process of natural linearisation.}

\bl\label{gradu}
 Let $\Phi\in L^2(0, \tau;  L^2(\O))$  be such that
$\nabla \Phi\in L^2(0, \tau;  L^\infty(\O,\R^d))$ and let $\C\in (L^2(\O))^{d\times d}$.  Then $\nabla\cdot \C\nabla \Phi \in L^2(0, \tau;   H^{-1}(\O)).$
\el
\bpf
Let $\C=(c_{ij})_{1\leq i,j\leq d}.$ First we show that $\C\nabla\Phi(\cdot,t)\in L^2(\O,\R^d)$ for $t\in [0,\tau].$ 

Since $\nabla \Phi(\cdot,t)\in L^\infty(\O,\R^d)$, for a.a $t\in [0,\tau]$, we have
$$
\Big\|{c_{ij}\dfrac{\partial\Phi(\cdot,t)}{\partial x_k}}\Big\|_{L^2}\leq \lVert{c_{ij}}\rVert_{L^2}\Big\|{\dfrac{\partial \Phi(\cdot,t)}{\partial x_k}}\Big\|_{L^\infty }\q \text{for all}\ 1\leq i,j,k\leq d\ \text{and a.a}\ t\in [0,\tau]. $$
Using Cauchy Schwarz inequality we have,
$$
\Big\|{\sum_{j=1}^dc_{ij}\dfrac{\partial\Phi(\cdot,t)}{\partial x_j} }\Big\|^2_{L^2}\leq \Big(\sum_{j=1}^d\lVert {c_{ij}}\rVert_{L^2}\Big\|{\dfrac{\partial\Phi(\cdot,t)}{\partial x_j} }\Big\|_{L^\infty}\Big)^2\leq\sum_{j=1}^d\lVert{c_{ij}}\rVert^2_{L^2} \sum_{j=1}^d\Big\|{\dfrac{\partial\Phi(\cdot,t)}{\partial x_j}}\Big\|^2_{L^\infty}$$
for a.a $t\in [0,\tau]$. 
Thus, we have
 $$\lVert{\C\nabla\Phi(\cdot,t)}\rVert^2_{L^2(\O,\R^d)}=\sum_{i=1}^d\Big\|{\sum_{j=1}^dc_{ij}\dfrac{\partial\Phi(\cdot,t)}{\partial x_j}}\Big\|_{L^2}^2\leq \sum_{i=1}^d\sum_{j=1}^d\lVert{c_{ij}}\rVert^2_{L^2} \sum_{j=1}^d\Big\|{\dfrac{\partial\Phi(\cdot,t)}{\partial x_j}}\Big\|^2_{L^\infty},$$
for a.a $t\in [0,\tau],$  and hence $$\lVert{\C\nabla \Phi(\cdot,t)}\rVert_{L^2(\O,\R^d)}\leq \lVert{\C}\rVert_{F}\lVert{\nabla \Phi(\cdot,t)}\rVert_{L^\infty(\O,\R^d)}.$$
Hence, $\C\nabla\Phi(\cdot,t)\in L^2(\O,\R^d)$ for a.a $t\in [0,\tau].$
Therefore, $\nabla \cdot \C\nabla \Phi(\cdot,t)\in H^{-1}(\O) $ for a.a $t\in [0,\tau].$
Now, using Lemma \ref{h-1bound}, for a.a $t\in [0, \tau]$, we have
\beqarray
\Big\|{\sum_{i=1}^d\sum_{j=1}^d}\dfrac{\partial}{\partial x_i}\Big(c_{ij}\dfrac{\partial\Phi(\cdot,t)}{\partial x_j}\Big)\Big\|_{H^{-1}}&\leq &  \sum_{i=1}^d\sum_{j=1}^d\Big\|{\dfrac{\partial}{\partial x_i}\Big(c_{ij}\dfrac{\partial\Phi(\cdot,t)}{\partial x_j}\Big)}\Big\|_{H^{-1}}\\
&\leq &\sum_{i=1}^d\sum_{j=1}^d\Big\lVert{c_{ij}\dfrac{\partial\Phi(\cdot,t)}{\partial x_j}}\Big\rVert_{L^2}\\
&\leq& \sum_{i=1}^d\sum_{j=1}^d\lVert{c_{ij}}\rVert_{L^2}\Big\|{\dfrac{\partial\Phi(\cdot,t)}{\partial x_j}}\Big\|_{L^\infty}\\
&\leq & d \lVert{\C}\rVert_{F}\lVert{\nabla \Phi(\cdot,t)}\rVert_{L^\infty(\O,\R^d)}.
\eeqarray
Thus, \beq\label{h-1estimate}
\lVert{\nabla\cdot \C\nabla\Phi(\cdot,t)}\rVert_{H^{-1}(\O)}\leq d \lVert{\C}\rVert_{F}\lVert{\nabla \Phi(\cdot,t)}\rVert_{L^\infty(\O,\R^d)}\q \text{for a.a}\ t\in[0,\tau].
\eeq
Since $\nabla \Phi\in L^2(0, \tau;  L^\infty(\O,\R^d)),$ we conclude that $\nabla\cdot \C\nabla\Phi \in L^2(0, \tau;  H^{-1}(\O)).$
\epf

Let $\A_0$ be as in (\ref{v0}). For $\C\in (L^2(\O))^{d\times d}$ and $w\in L^2(0,\tau;L^2(\O)) $, we  consider the PDE
\beq\label{pde4}
v_t-\nabla\cdot \A_0\nabla v=\nabla\cdot \C\nabla w\q \text{in}\ Q,
\eeq
along with the  conditions
\beq\label{bdry4}
v=0\q \text{on} \ \partial Q,
\eeq
\beq\label{initial4}
v(x,0)=0\q\text{in}\ \O.
\eeq
Note that to talk about the existence of a weak solution of \eqref{pde4}-\eqref{initial4}, it is enough to make sure that the R.H.S of \eqref{pde4}, that is, $\nabla\cdot \C\nabla w$ belongs to $ L^2(0, \tau;  H^{-1}(\O)).$ In order to satisfy this condition, Lemma \ref{gradu} suggests that, it is enough to assume that $\nabla w\in L^2(0, \tau;  L^\infty(\O,\R^d))$. Thus, under this assumption, by Theorem \ref{dirichletestimate}, for each $\C\in (L^2(\O))^{d\times d}$, there exist a unique weak solution $v\in \X$, of \eqref{pde4}-\eqref{initial4}, where $\X$ is defined as in (\ref{space-new}).

{In view of the discussion in the above paragraph,  for each  $w\in L^2(0,\tau;L^2(\O))$ with  $\nabla w\in L^2(0, \tau;  L^\infty(\O,\R^d))$, we consider the map
 $\T_w:(L^2(\O))^{d\times d} \to \X$  defined by
\beq\label{tdef}\T_w \C=v,\q \C\in (L^2(\O))^{d\times d},\eeq
where  $v\in \X$ is the unique weak solution of \eqref{pde4}-\eqref{initial4}. 
Clearly,  the map $\T_w$ is a linear operator.  We now have the following result.
}

\bt\label{tprop1}
Let ${w\in L^2(0,\tau;L^2(\O))}$ be such that  {$0\neq\nabla w\in L^2(0, \tau;  L^\infty(\O,\R^d))$}.  {Then  the linear operator $\T_w:(L^2(\O))^{d\times d} \to \X$,
defined in \eqref{tdef},}  is a bounded linear operator of infinite rank with $$\lVert{\T_w}\rVert\leq d\sqrt{C_0}\Big(\int_0^\tau\lVert{\nabla w(\cdot,t)}\rVert^2_{L^\infty}dt\Big)^{1/2},$$
where  $C_0$ is  the constant as in Theorem \ref{dirichletestimate}.
\et

\bpf
We first show that $\T_w$ is a bounded operator. Let $\C\in (L^2(\O))^{d\times d}$ and $\T_w\C=v.$ Then, {$v$ is the unique weak solution of \eqref{pde4}-\eqref{initial4}.}  Hence using the estimate \eqref{dirichletestimatev} of Theorem \ref{dirichletestimate}, we have
\beqarray
\lVert{\T_w\C}\rVert=\lVert{v}\rVert_{\X}\leq C_0 \lVert{\nabla\cdot \C\nabla u}\rVert_{L^2(0, \tau;  H^{-1}(\O))}.
\eeqarray
Now, using the estimate given in \eqref{h-1estimate}, we have
\beqarray
\lVert{\nabla\cdot \C\nabla w}\rVert^2_{L^2(0, \tau;  H^{-1}(\O))}&= &\int_0^\tau\lVert{\nabla\cdot \C\nabla w(\cdot,t)}\rVert^2_{H^{-1}(\O)}dt\\
&\leq & d^2\lVert{\C}\rVert^2_{F}\int_0^\tau\lVert{\nabla w(\cdot,t)}\rVert^2_{L^\infty(\O)}dt.
\eeqarray
Thus, $$\lVert{\T_w\C}\rVert\leq d\sqrt{C_0}\big(\int_0^\tau\lVert{\nabla w(\cdot,t)}\rVert^2_{L^\infty}dt\big)^{1/2}\lVert{\C}\rVert_{F}. $$ This shows that $\T_w$ is a bounded linear operator with $\lVert{\T_w}\rVert\leq d\sqrt{C_0}\big(\int_0^\tau\lVert{\nabla w(\cdot,t)}\rVert^2_{L^\infty}dt\big)^{1/2}. $

We now show that $\T_w$ {is of}  infinite rank. Since $\nabla w(\cdot,t)\neq 0$ for almost all $t\in [0,\tau],$ without loss of generality we assume that $\dfrac{\partial w}{\partial x_1}\neq 0$ for almost all $t\in [0,\tau]$. Let  $(\varphi_n)$ be a sequence of elements in $C_c^\infty(\O)$ such that ${\rm supp}\ \varphi_m\cap\ {\rm supp}\ \varphi_n=\phi$ for $m\neq n.$
Let $\C_n$ be the $d\times d$ matrix given by,
$$\C_n=\begin{pmatrix}
\varphi_n &O_{(d-1)\times(d-1)}\\
O_{1\times (d-1)}&0
\end{pmatrix},
$$
where $O_{k\times l}$ denotes the zero matrix of respective order.
Then, clearly $(\C_n)$ is a sequence of linearly independent elements in $(L^2(\O))^{d\times d}.$ Let $\T_w\C_n=v_n.$ We claim that $\{v_n:n\in \N\}$ is linearly independent. {Assume, for a moment, that $\{v_n:n\in \N\}$ is linearly dependent.} Then without loss of generality we assume that, for $m\in \N,$ let  $v_1=\sum_{i=2}^m\b_i v_i,$ where $\b_i$'s are constants not all zero.  Using the definition of $\T_w$, \eqref{pde4} gives
\beqarray
\nabla\cdot \C_1\nabla w=(v_1)_t-\nabla \cdot \A_0\nabla v_1=\sum_{i=2}^m\b_i(v_i)_t-\sum_{i=2}^m\b_i\nabla \cdot \A_0\nabla v_i
=\sum_{i=2}^m\b_i\nabla\cdot \C_i\nabla w.
\eeqarray
Thus, $\ds \nabla\cdot (\C_1-\sum_{i=2}^m\b_i \C_i)\nabla w=0.$ Hence, we have
\beq\label{lid}
\int_\O(\C_1-\sum_{i=2}^m\b_i\C_i)\nabla w\cdot\nabla\f\,dx=0\q\text{for all }\,\, \f\in H^1_0(\O)\,\, \text{and for a.a.}\,\, t\in (0,\tau). 
\eeq
{ Let $\f\in H^1_0(\O)$ be such that $\frac{\partial\f}{\partial x_1}=\varphi_1\dfrac{\partial w}{\partial x_1}-\sum_{i=2}^m\b_i \varphi_i\dfrac{\partial w}{\partial x_1}$.
Then from \eqref{lid}, we obtain
$$\int_\O(\f_1-\sum_{i=2}^m\b_i\f_i)^2\Big(\frac{\partial w}{\partial x_1}\Big)^2\,dx=0\q\text{for a.a.}\,\, t\in (0,\tau).$$ Since $\frac{\partial w}{\partial x_1}\neq 0$ a.e. in $(0,\tau),$ we have $\f_1=\sum_{i=2}^m\b_i\f_i$  a.e. in $\O$.  This leads to a contradiction, since $\f_n$'s are linearly independent.}
Therefore, $\{v_n:n\in \N\}$ is a linearly independent set. This shows that $\T_w$ is of infinite rank.
\epf

Our next theorem demonstrates  one more property of the operator $\T_w.$

\bt\label{tcmpt}
Let $w$ be as assumed in Theorem \ref{tprop1} and $\T_w:(L^2(\O))^{d\times d}\to \X$ be as defined in \eqref{tdef}. Then $\T_w:(L^2(\O))^{d\times d}\to L^2(0, \tau;  L^2(\O))$ is a compact linear operator.
\et

\bpf
By {Theorem \ref{tprop1},  $\T_w:(L^2(\O))^{d\times d}\to \X$}  is a bounded linear operator. Also, by Remark \ref{compactembedding}, $\X$ is compactly embedded in $L^2(0, \tau;  L^2(\O)).$ Therefore $\T_w:(L^2(\O))^{d\times d}\to L^2(0, \tau;  L^2(\O))$ is a composition of a bounded linear operator and a compact linear operator. Hence $\T_w:(L^2(\O))^{d\times d}\to L^2(0, \tau;  L^2(\O))$ is a compact linear operator.
\epf

 We are now in a position to represent  our non linear inverse problem in an alternative form as 
 \beq\label{lininvopeq}
\T_u\B=v_0-u,
\eeq
where  $u$ is  as in the Assumption \ref{existence},  the linear operator $\T_u:(L^2(\O))^{d\times d}\to L^2(0,\tau;L^2(\O))$ is as defined in \eqref{tdef},   $v_0$ is as in \eqref{v0}, and $\B =\A-\A_0$ is a solution of \eqref{lininvopeq},  where $\A_0$ is as considered in Theorem \ref{dirichletestimate}.

{Since the operator $\T_u$ depends on the data $u$ and it is  a compact linear operator of infinite rank (cf. Theorem \ref{tprop1} and Theorem \ref{tcmpt}), the operator equation \eqref{lininvopeq} is  ill-posed, that is, small perturbations in the data $u$ may lead to a large deviation in the solution. But, in practical application, noise in the data is inevitable. So, a regularization method is necessary for obtaining stable approximate solution, and we employ the method of Tikhonov regularization for the same. For more details about regularization theory, one may look into the books \cite{neuber,  nairopeq}.}
}

\section{Regularization and error analysis}\label{reg}
As mentioned in the previous section,  $\B:=\A-\A_0$ satisfies the equation \eqref{lininvopeq}.   Thus, we have
\beq\label{exacteqn}
\T_u\B=v_0-u.
\eeq
Now, suppose that $u$ is available with some noise, say we have $\tilde u$ in place of $u$.  
\noindent
{\bea\label{noise}  For  $\d>0$, let  $\tilde{u}\in {L^2(0,\tau;L^2(\O))}$ be the observed data corresponding to $u,$ satisfying $$\nabla \tilde{u}\in L^2(0, \tau;  L^\infty(\O,\R^d))$$
along with
\beq\label{uestimate}
\lVert{\tilde{u}-u}\rVert_{L^2(0,\tau;L^2(\O))}+ \Big(\int_0^\tau \lVert{\nabla u(\cdot,t)-\nabla \tilde{u}(\cdot,t)}\rVert^2_{L^\infty}dt\Big)^{1/2}\leq \d\ .
\eeq
\ena

\brem\label{smoothingremark}
Note that in the above assumption we have imposed a regularity condition on the observed data $\tilde{u}$, namely $\nabla \tilde{u}\in L^2(0,\tau;L^\infty(\O,\R^d)).$ But, $\tilde{u}$ being an observed data, may not satisfy that regularity condition. In order to accommodate the noisy data $\tilde{u}$ only  with the condition that it belongs to ${L^2(0,\tau;L^2(\O))}$ and
$$\lVert{\tilde{u}-u}\rVert_{L^2(0,\tau;L^2(\O))}\leq \d$$
we consider a a \textit{smoothed version} of $\tilde{u}$, say $\tilde{z}$, which will satisfy all the conditions in the above assumption and then carry out the analysis with $\tilde{z}$ in place of $\tilde{u}$. In Section \ref{smoothingnoisy}, we have given the procedure for obtaining the smoothed version of noisy data. 
\erem

{\bf Notation:} In rest of the paper whenever needed, we shall use the notation $\|u-\tilde{u}\|$ to denote $\|u-\tilde{u}\|_{L^2(0,\tau;L^2(\O))}$.

{
Corresponding to the noisy data  $\tilde{u}$ satisfying the  Assumption \ref{noise}, we would like to obtain a  regularized approximation  of  $\B$. For this purpose we shall consider the Tikhonov regularization (cf. \cite{neuber, nairopeq}).  That is, for each $\a>0$,  the candidate for the regularized approximation  of  $\B:=\A-\A_0$ is the unique solution of the equation 
\beq\label{tikhpertbd1}
(\T_{\tilde{u}}^* \T_{\tilde{u}}+\a I)\tilde{\B}_\a=\T_{\tilde{u}}^*(v_0-\tilde{u}).
\eeq
Then we take  the regularized approximation of $\A$ as 
$$ \tilde{\A}_\a =  \A_0+\tilde{\B}_\a,\q \a>0.$$
 Let  $\B_\a$ be the unique solution of
\beq\label{tikhexact2}
(\T_u^*\T_u+\a I)\B_\a= \T_u^*(v_0-u).
\eeq
Also, using \eqref{exacteqn}, we have
\beq\label{tikhpertbd2}
(\T_u^*\T_u+\a I)\B=\T_u^*(v_0-u)+\a \B.
\eeq 
Then by Theorem \ref{tprop1} and using the estimate given in \eqref{uestimate}, we have$$\lVert{\T_u\B-\T_{\tilde{u}}\B}\rVert\leq d\sqrt{C_0}\Big(\int_0^\tau \lVert{\nabla u(\cdot,t)-\nabla \tilde{u}(\cdot,t)}\rVert^2_{L^\infty}dt\Big)^{1/2}\lVert{\B}\rVert_{F}$$which implies that
\beq\label{estimateoft}
\lVert{\T_u-\T_{\tilde{u}}}\rVert\leq d\sqrt{C_0}\Big(\int_0^\tau \lVert{\nabla u(\cdot,t)-\nabla \tilde{u}(\cdot,t)}\rVert^2_{L^\infty}dt\Big)^{1/2}\leq d\sqrt{C_0}\d,
\eeq
where  $C_0$ be as in \eqref{dirichletestimatev}. 
}
The equations \eqref{tikhpertbd1} and \eqref{tikhpertbd2} leads to
\begin{eqnarray}
\tilde{\B}_\a-\B_\a &=&(\T_{\tilde{u}}^*\T_{\tilde{u}}+\a I)^{-1}\T_{\tilde{u}}^*(v_0-\tilde{u})-(\T_u^*\T_u+\a I)^{-1}\T_u^*(v_0-u)\nonumber\\
&=&[(\T_{\tilde{u}}^*\T_{\tilde{u}}+\a I)^{-1}\T_{\tilde{u}}^*-(\T_u^*\T_u+\a I)^{-1}\T_u^*](v_0-u)\label{subham}\\
&&+(\T_{\tilde{u}}^*\T_{\tilde{u}}+\a I)^{-1}\T_{\tilde{u}}^*(u-\tilde{u}).\nonumber
\end{eqnarray}
Now, using \eqref{exacteqn}, we have
$$(\T_{\tilde{u}}^*\T_{\tilde{u}}+\a I)^{-1}\T_{\tilde{u}}^* (v_0-u) = (\T_{\tilde{u}}^*\T_{\tilde{u}}+\a I)^{-1}\T_{\tilde{u}}^*\T_u\B
$$
$$(\T_u^*\T_u+\a I)^{-1}\T_u^*(v_0-u) = (\T_u^*\T_u+\a I)^{-1}\T_u^*\T_u\B
$$
Also, we see that
\beqarray
(\T_{\tilde{u}}^*\T_{\tilde{u}}+\a I)^{-1}\T_{\tilde{u}}^*\T_u - (\T_u^*\T_u+\a I)^{-1}\T_u^*\T_u
&=&  (\T_{\tilde{u}}^*\T_{\tilde{u}}+\a I)^{-1}\T_{\tilde{u}}^*(\T_u-\T_{\tilde{u}})\T_u^*\T_u (\T_u^*\T_u+\a I)^{-1}\\
&& +\a (\T_{\tilde{u}}^*\T_{\tilde{u}}+\a I)^{-1}(\T_{\tilde{u}}^*-\T_u^*)(\T_u\T_u^*+\a I)^{-1}\T_u.
\eeqarray
Hence, (\ref{subham}) takes the form
\begin{eqnarray*}
\tilde{\B}_\a-\B_\a &=&
(\T_{\tilde{u}}^*\T_{\tilde{u}}+\a I)^{-1}\T_{\tilde{u}}^*(\T_u-\T_{\tilde{u}})\T_u^*\T_u (\T_u^*\T_u+\a I)^{-1}\B\\
&& +\a (\T_{\tilde{u}}^*\T_{\tilde{u}}+\a I)^{-1}(\T_{\tilde{u}}^*-\T_u^*)(\T_u\T_u^*+\a I)^{-1}\T_u\B\\
&&+(\T_{\tilde{u}}^*\T_{\tilde{u}}+\a I)^{-1}\T_{\tilde{u}}^*(u-\tilde{u}).
\eeqarray
Now, we recall the following estimates (cf. \cite{nairopeq}):
\beqarray
\lVert{(\T_{\tilde{u}}^*\T_{\tilde{u}}+\a I)^{-1}\T_{\tilde{u}}^*}\rVert\leq \dfrac{1}{2\sqrt{\a}},&&\q \lVert{(\T_{\tilde{u}}^*\T_{\tilde{u}}+\a I)^{-1}}\rVert\leq \dfrac{1}{\a}\\
\lVert{\T_u^*\T_u(\T_u^*\T_u+\a I)^{-1}}\rVert\leq 1,&&\q \lVert{(\T_u\T_u^*+\a I)^{-1}\T_u}\rVert\leq \dfrac{1}{2\sqrt{\a}}.
\eeqarray
Using these estimates, we obtain $$\lVert{[(\T_u^*\T_u+\a I)^{-1}\T_u^*-(\T_{\tilde{u}}^*\T_{\tilde{u}}+\a I)^{-1}\T_{\tilde{u}}^*](v_0-u)}\rVert\leq \dfrac{\lVert{\T_u-\T_{\tilde{u}}}\rVert}{\sqrt{\a}}\lVert{\B}\rVert$$
and $$\lVert{(\T_{\tilde{u}}^*\T_{\tilde{u}}+\a I)^{-1}\T_{\tilde{u}}^*(\tilde{u}-u)}\rVert\leq \dfrac{\lVert{\tilde{u}-u}\rVert}{2\sqrt{\a}}.$$
Thus, from (\ref{subham}), we have $$\lVert{\B_\a-\tilde{\B}_\a}\rVert\leq \dfrac{\lVert{\T_u-\T_{\tilde{u}}}\rVert}{\sqrt{\a}}\lVert{\B}\rVert+\dfrac{\lVert{\tilde{u}-u}\rVert}{2\sqrt{\a}}. $$
Therefore, by \eqref{uestimate}, we have $$\lVert{\B_\a-\tilde{\B}_\a}\rVert\leq \dfrac{d\sqrt{C_0}\d}{\sqrt{\a}}\lVert{\B}\rVert+\dfrac{\d}{2\sqrt{\a}}=(d\sqrt{C_0}\|\B\|+1/2)\dfrac{\d}{\sqrt{\a}}.$$
It is well known (cf. \cite{nairopeq}) from the theory of Tikhonov regularization that
$$\lVert{\B-\B_\a}\rVert\to 0\q\h{as}\q \a\to 0.$$
Thus, {since  $\A-\tilde{\A}_\a = \B-\tilde{\B}_\a$ for  $\a>0$, }
we have proved the following result.

\bt\label{error}
Let $\A_0$ be as in Theorem \ref{dirichletestimate} and $\B$ be as given in \eqref{exacteqn}. Let $\tilde{\B}_\a$ be the unique solution of \eqref{tikhpertbd1} and $\B_\a$ be as given in \eqref{tikhexact2}. Let $\d$ and $\tilde{u}$ be as given in \eqref{uestimate}. Then we have,
$$\lVert{\A-\tilde{\A}_\a}\rVert\leq \lVert{\B-\B_\a}\rVert+(d\sqrt{C_0}\|\B\|+1/2)\dfrac{\d}{\sqrt{\a}},$$
where  $\lVert{\B-\B_{\a}}\rVert\to 0$ as $\a\to 0.$ Further, 
choose the parameter $\a:=\a_\d$ depending on $\d$ in such a way that $\a_\d\to 0$ and
$$\dfrac{d\sqrt{C_0}\,\d}{\sqrt{\a_\d}}\lVert{\B}\rVert+\dfrac{\d}{2\sqrt{\a_\d}}\to 0\q\h{as}\q \d\to 0,$$
then we have $\lVert{\A-\tilde{\A}_{\a_\d}}\rVert\to 0$ as $\d\to 0.$
\et

In the next subsection, we determine an explicit  representation of  f$\T_w^*$, the adjoint of $\T_w$.

\subsection{Explicit representation of the adjoint}
Let  $\T_w:(L^2(\O))^{d\times d}\to L^2(0, \tau;  L^2(\O))$  be as in Theorem \ref{tcmpt}.  
Let $\C\in (L^2(\O))^{d\times d}$ and $v=\T_w\C.$ Then,{ 
 in view of  (\ref{bdry4}) - (\ref{tdef}), $v\in L^2(0,\tau;H^1_0(\O))$  and we have}
\beq\label{weak-adj-v-2}
\<v_t,\f\>+\int_\O\A_0\nabla v\cdot \nabla \f\,dx=-\int_\O\C\,\nabla w\cdot \nabla \f\,dx\q\forall\, \f\in L^2(0, \t; H^1_0(\O))\h{ and for a.a. }\, t\in [0,\tau],
\eeq
{where $v_t\in L^2(0, \t; H^{-1}(\O))$, and $\<\cdot, \cdot\>$ in the above is the duality pairing corresponds to   $H^{-1}(\O)$ and $H^{1}_0(\O)$.}
Now, let   $\phi\in L^2(0,\tau;L^2(\O))$ and consider the PDE
\beq\label{adj-pde-2}
\begin{cases}
z_t+\nabla\cdot\A_0\nabla z=\phi \q&\text{in}\,\,Q \\
z=0\q&\text{on}\,\,\partial Q\\
z(\cdot,\tau)=0\q& \text{in}\,\,\O.
\end{cases}
\eeq
By  reversing the time direction for the above PDE, and by Theorem \ref{dirichletestimate}, {we know that there exists a unique $z\in L^2(0,\tau;H^1_0(\O))$ with $z_t\in L^2(0,\tau;H^{-1}(\O))$ such that  $z(\cdot,\tau)=0$}  and
\beq\label{weak-adj-z-2}
\<z_t,\f\> - \int_\O\A_0\nabla z\cdot\nabla\f\,dx=\int_\O\phi\f\,dx\q\text{for a.a.}\,\, t\in [0,\tau]\,\, \text{and}\,\,\f\in H^1_0(\O).
\eeq 
Since $v(\cdot,t)\in H^1_0(\O)$ for a.a. $t\in [0,\tau],$ from \eqref{weak-adj-z-2}, we obtain
$$\<z_t,v\> -\int_\O\A_0\nabla z\cdot\nabla v\,dx=\int_\O\phi v\,dx\q\text{for a.a.}\,\, t\in [0,\tau].
$$
Now, integrating  with respect to $t$, and using the fact $v(\cdot,0)=0=z(\cdot,\tau),$ we obtain
\beq\label{adj-var1-2}
-\int_0^\tau\<z,v_t\>dt-\int_0^\tau\int_\O\A_0\nabla z\cdot \nabla v\,dxdt=\int_0^\tau\int_{\O}\phi v\,dxdt=\<\T_w\C,\phi\>_{L^2(0,\tau;L^2(\O))}.
\eeq
Since $z\in L^2(0,\tau;H^1_0(\O)),$ from \eqref{weak-adj-v-2}, we obtain
$$\<v_t,z\> +\int_\O\A_0\nabla v\cdot \nabla z\,dx=-\int_\O\C\nabla w\cdot \nabla z\,dx\q\text{for a.a.}\,\, t\in [0,\tau],$$and hence
\beq\label{adj-var2-2}
\int_0^\tau\<v_t,z\>dt+\int_0^\tau\int_\O\A_0\nabla v\cdot \nabla z\,dxdt=-\int_0^\tau\int_\O\C\nabla w\cdot \nabla z\,dxdt.
\eeq
Therefore, from \eqref{adj-var1-2} and \eqref{adj-var2-2}, using the fact that $\A_0$ is symmetric, we obtain
\beq\label{adj-final-2}
\int_0^\tau\int_{\O}\phi v\,dxdt=\int_0^\tau\int_\O\C\nabla w\cdot \nabla z\,dxdt.
\eeq
Now, 
\beqarray
\int_0^\tau\int_\O\C\nabla w\cdot \nabla z\,dxdt&=&\int_0^\tau\int_\O\<\C\nabla w,\nabla z\>_{\R^d}\,dxdt=\int_0^\tau\int_\O(\C\nabla w)^T\nabla z\,dxdt\\
&=& \int_0^\tau\int_\O (\nabla w)^T\C^T\nabla z\,dxdt=\int_0^\tau\int_\O (\nabla w)^T\big((\nabla z)^T\C\big)^T\,dxdt\\
&=& \int_0^\tau\int_\O \<(\nabla z)^T\C,(\nabla w)^T\>_{\R^{1\times d}}\,dxdt=\int_\O \<\C, \int_0^\tau\nabla z(\nabla w)^Tdt\>_{\R^{d\times d}}\,dx\\
&=& \Big\<\C,\int_0^\tau\nabla z(\nabla w)^Tdt\Big\>_{(L^2(\O))^{d\times d}}.
\eeqarray 
Thus, from \eqref{adj-final-2}, we have
$$\<T_w\C,\phi\>_{L^2(0,\tau;L^2(\O))}=\int_0^\tau\int_{\O}\phi v\,dxdt=\Big\<\C,\int_0^\tau \nabla z(\nabla w)^Tdt\Big\>.$$
Therefore, we have proved the following result.

\bt\label{adj-rep-2}
Let $w$ be as in Theorem \ref{tprop1} and $\T_w:(L^2(\O))^{d\times d}\to L^2(0,\tau;L^2(\O))$ be as  in Theorem \ref{tcmpt}. Then the operator $\T_w^*,$ the adjoint of $\T_w$,  is given by
$$\T_w^*\phi=\int_0^\tau\nabla z(\nabla w)^Tdt,\q\phi\in L^2(0,\tau;L^2(\O)),$$
where $z\in L^2(0,\tau;H^1_0(\O))$ is the unique weak solution of the PDE:
$$z_t+\nabla\cdot\A_0\nabla z=\phi\,\,\text{in}\,\,\O\times (0,\tau)$$
along with the conditions
$$z=0\,\,\text{on}\,\,\partial\O\times (0,\tau),\q z(\cdot,\tau)=0\,\,\text{in}\,\,\O.$$
\et
}
%%%%%%%
%%%%%%%
\section{Finite dimensional setting}\label{finit}
So far we have theoretically obtained stable approximation for $\A$, {namely $\tilde{\A}_\a:= \A_0+\tilde{\B}_\a$,} and we have seen our underlying spaces were all infinite dimensional. But, {in the context of applications  one has to} work with finite dimensional spaces. From that point of view it becomes necessary to realize our analysis in finite dimensional setting. For this purpose,  we employ Galerkin method to obtain finite dimensional approximations of $\tilde{\B}_\a,$ which in turn will give approximations for $\A$ with suitable choice of the regularization parameter $\a$.

Let $(X_n)$ be a sequence of finite dimensional subspaces of $L^2(\O)$ such that
\ben
\i[(i)] $X_n\subseteq X_{n+1}$ for all $n\in \N$ and
\i[(ii)] $\bigcup_{n=1}^\infty X_n$ is dense in $L^2(\O).$
\een
For each $n\in \N,$ let $P_n:L^2(\O)\to L^2(\O)$ be an orthogonal projection onto $X_n.$ Then, with the assumptions (i) and (ii) on $X_n,$ we have $$P_n\varphi\to \varphi\q\h{as}\q  n\to \infty$$ 
for all $\varphi\in L^2(\O).$
For each $n\in \N,$ we define $\P_n:(L^2(\O))^{d\times d}\to (L^2(\O))^{d\times d}$ by
\beq\label{matrixprojection}
\P_n\C=(P_nc_{ij})_{1\leq i,j\leq d}
\eeq
for all $\C=(c_{ij})_{1\leq i,j\leq d}\in (L^2(\O))^{d\times d}.$ 

We now have the following result.
{
\bt\label{matrixprojectionconv}
For each $n\in \N$, let  $\P_n:(L^2(\O))^{d\times d}\to (L^2(\O))^{d\times d}$  be as defined in \eqref{matrixprojection}. Then $\P_n$ 
is an orthogonal projection and for every $\C\in (L^2(\O))^{d\times d},$ 
$$\lim_{n\to \infty}\lVert{\P_n\C-\C}\rVert\to 0.$$
\et 
\bpf
From the definition of $\P_n,$ it follows that $\P_n$ is a projection. Also, $X_n\subseteq X_{n+1}$ for all $n\in \N$ ensures that $R(\P_n)\subseteq R(\P_{n+1})$ for all $n\in \N,$ where $R(\P_n)$ denotes the range space of $\P_n.$
Since $P_n:L^2(\O)\to L^2(\O)$ is an orthogonal projection,
for any $\C,\, \D\in (L^2(\O))^{d\times d}$, we have
$$\Big\langle\C,\P_n\D\Big\rangle=\sum_{i,j=1}^d\langle c_{ij}, P_nd_{ij}\rangle_{L^2}=\sum_{i,j=1}^d\langle P_nc_{ij},d_{ij}\rangle _{L^2}= \Big\langle\P_n\C,\D\Big\rangle.$$
Hence, $\P_n:(L^2(\O))^{d\times d}\to (L^2(\O))^{d\times d}$ is an orthogonal projection.
Also, we have 
$$\lVert{\P_n\C-\C}\rVert^2 =\sum_{i,j=1}^d\lVert{P_nc_{ij}-c_{ij}}\rVert_{L^2}^2.$$ Since $P_n$ converges pointwise to the identity  in $L^2(\O),$ we have $$\lim_{n\to \infty}\lVert{P_nc_{ij}-c_{ij}}\rVert_{L^2}\to 0\q \text{for all}\ 1\leq i,j\leq d.$$Therefore, we have $\ds \lim_{n\to \infty}\lVert{\P_n\C-\C}\rVert=0$ for every $\C\in (L^2(\O))^{d\times d}.$
\epf
}

{
\bcor\label{cpt-appr} 
For each $n\in \N$, let  $\P_n:(L^2(\O))^{d\times d}\to (L^2(\O))^{d\times d}$  be as defined in \eqref{matrixprojection}. Then  for every $w\in L^2(0, \t:, L^2(\O))$ with $0\neq\nabla w\in L^2(0, \tau;  L^\infty(\O,\R^d))$, 
$$\|\T_w - \T_w \P_n\|\to 0 \q\h{as}\q n\to \infty.$$ 
\ecor

\bpf
We know, by Theorem \ref{tcmpt}, that $\T_w$ is a compact operator and by Theorem \ref{matrixprojectionconv} that $\P_n$ is an orthogonal projection for each $n\in \N$ satisfying  $\lim_{n\to \infty}\lVert{\P_n\C-\C}\rVert\to 0. $ 
Since  $\T_w^*$ is also a compact operator, by a standard result in Functional analysis, as a consequence of uniform boundedness principle (see, e.g., \cite{nair-fa},  Corollary 6.6),  we have    
$$\|(I-\P_n) \T_w^*\| \to 0 \q\h{as}\q n\to \infty.$$
Since   
$\|\T_w - \T_w \P_n\| = \|(I-\P_n) \T_w^*\|$, the conclusion in the corollary follows.   
\epf 
}

Next we observe that \eqref{tikhpertbd1} holds iff
\beq\label{tikhpert1ipform}
\Big\langle (\T_{\tilde{u}}^* \T_{\tilde{u}}+\a I)\tilde{\B}_\a,\C\Big\rangle=\Big\langle\T_{\tilde{u}}^*(v_0-\tilde{u}),\C\Big\rangle\q \text{for all}\ \C\in (L^2(\O))^{d\times d}.
\eeq
In order to obtain finite dimensional approximations of $\tilde{\B}_\a,$ in \eqref{tikhpert1ipform}, we vary $\C\in R(\P_n)$ for each $n\in \N.$ In other words, for each $n\in \N$ and $\a>0,$ we look for $\tilde{\B}_{\a, n}\in R(\P_n)$ such that
\beqarray
\Big\langle (\T_{\tilde{u}}^* \T_{\tilde{u}}+\a I)\tilde{\B}_{\a, n},\C\Big\rangle=\Big\langle\T_{\tilde{u}}^*(v_0-\tilde{u}),\C\Big\rangle\q \text{for all}\ \C\in R(\P_n).
\eeqarray
Equivalently, we look for $\tilde{\B}_{\a, n}\in R(\P_n)$ such that
\beq\label{fintetikhpert1ipform}
 \Big\langle (\P_n\T_{\tilde{u}}^* \T_{\tilde{u}}\P_n+\a I)\tilde{\B}_{\a, n},\C\Big\rangle=\Big\langle\P_n\T_{\tilde{u}}^*(v_0-\tilde{u}),\C\Big\rangle\q \text{for all}\ \C\in (L^2(\O))^{d\times d}.
\eeq
Equivalently,
\beq\label{finitetikhpert1opeq}
(\P_n\T_{\tilde{u}}^* \T_{\tilde{u}}\P_n+\a I)\tilde{\B}_{\a, n}=\P_n\T_{\tilde{u}}^*(v_0-\tilde{u}).
\eeq
For each $n\in \N,$ let dim$(X_n)=n$. Then, it can be seen that dim$(R(\P_n))=nd^2.$ Indeed, if $\{\varphi_1,...,\varphi_n\}$ is a basis of $X_n$, then it can be seen that
$\{E^{l^{ij}}_{ij}: 1\leq i,j\leq d,\ 1\leq l^{ij}\leq n \}$ forms a basis for $R(\P_n),$ where
$$
(E^{l^{ij}}_{ij})_{pq}= \left\{\begin{array}{ll}
                                       \varphi_{l^{ij}}
&\h{if}\ p=i, q=j\\
                                       0
&\h{if}\ p\neq i\ \text{or}\ q\neq j. \end{array}\right.
$$
Also, it can be seen that, if $\{\varphi_1,...,\varphi_n\}$ is an orthonormal basis of $X_n$, then $\{E^{l^{ij}}_{ij}: 1\leq i,j\leq d,\ 1\leq l^{ij}\leq n\}$ too becomes an orthonormal basis for $R(\P_n).$

Now, let $\{\varphi_1,...,\varphi_n\}$ be an orthonormal basis of $X_n$, and let $\{E^{l^{ij}}_{ij}: 1\leq i,j\leq d,\ 1\leq l^{ij}\leq n\}$, which is  an orthonormal basis of $R(\P_n).$ Let us write the solution of the well-posed equation (\ref{finitetikhpert1opeq}) as
\beq\label{basisrep}
\tilde{\B}_{\a, n}=\sum_{1\leq i,j\leq d}\ \sum_{l^{ij}=1}^n c^{l^{ij}}_{ij}E^{l^{ij}}_{ij}
\eeq
for some constants $c^{l^{ij}}_{ij},\ 1\leq i,j\leq d,\ 1\leq l^{ij}\leq n.$
Then, from \eqref{finitetikhpert1opeq}, we have
$$\Big\langle(\P_n\T_{\tilde{u}}^* \T_{\tilde{u}}\P_n+\a I)\tilde{\B}_{\a, n},E^{l^{ij}}_{ij}\Big\rangle =\Big\langle\P_n\T_{\tilde{u}}^*(v_0-\tilde{u}), E^{l^{ij}}_{ij}\Big\rangle\q \text{for all}\  1\leq i,j\leq d,\ 1\leq l^{ij}\leq n. $$
Using \eqref{basisrep}, we have
\beq\label{sumform}
\sum_{1\leq p,q \leq d}\ \sum_{k^{pq}=1}^nc^{k^{pq}}_{pq}\langle \T_{\tilde{u}}E^{k^{pq}}_{pq}, \T_{\tilde{u}}E^{l^{ij}}_{ij}\rangle + \a \sum_{1\leq p,q\leq d}\ \sum_{k^{pq}=1}^n c^{k^{pq}}_{pq}\langle E^{k^{pq}}_{pq}, E^{l^{ij}}_{ij}\rangle =\langle (v_0-\tilde{u}),\T_{\tilde{u}}E^{l^{ij}}_{ij}\rangle
\eeq
for all $1\leq i,j\leq d,\ 1\leq l^{ij}\leq n.$ The equation  \eqref{sumform} can be written in the matrix form as
\beq\label{reducedmatrix}
U \vec{c}+\a D \vec{c}=\vec{b},
\eeq
where $U, D$ are matrices and $\vec{c},\ \vec{b}$ are column vectors  given by
$$U={\bf \Big[}\langle \T_{\tilde{u}}E^{k^{pq}}_{pq},\,  \T_{\tilde{u}}E^{l^{ij}}_{ij}\rangle  {\bf \Big]},\q \ D= {\bf\Big[}\langle E^{k^{pq}}_{pq},\,   E^{l^{ij}}_{ij}\rangle {\bf \Big]},\q \vec{c}=  {\bf\Big[}c^{k^{pq}}_{pq}{\bf\Big]},\,  \vec{b}={\bf\Big[}\langle (v_0-\tilde{u}),\,  \T_{\tilde{u}}E^{l^{ij}}_{ij}\rangle {\bf\Big]} $$
for $1\leq p,q\leq d,\ 1\leq k^{pq}\leq n\q \text{and}\q 1\leq i,j\leq d,\ 1\leq l^{ij}\leq n.$

Thus, in order to obtain $\tilde{\B}_{\a, n}$, we need to solve the matrix equation \eqref{reducedmatrix}. Since $\P_n\T_{\tilde{u}}^*\T_{\tilde{u}}\P_n$ is a bounded {positive} self adjoint operator, for each $\a>0,$ there exists a unique solution for the equation \eqref{finitetikhpert1opeq}, in other words $\tilde{\B}_{\a, n}$ is determined uniquely and hence the matrix equation \eqref{reducedmatrix}, also has a unique solution. Thus, $\vec{c}$ is determined uniquely and hence $\tilde{\B}_{\a, n}$ is obtained explicitly using \eqref{basisrep}.

Recall that our goal is to obtain finite dimensional approximations to $\tilde{\B}_\a.$ So far we have obtained $\tilde{\B}_{\a, n}$, which is obtained by solving some matrix equation. But, we do not know whether these $\tilde{\B}_{\a, n}$ will approximate $\tilde{\B}_\a$ in some sense or not. 
%To settle this question about approximation, we have the following result. 
{Now, in view of Corollary \ref{cpt-appr}, the following theorem shows that  $\tilde{\B}_{\a, n}$ is an approximation of $\B_\a$ for each $\a$ if $n$ is large enough and $\d$ is small enough. }

\bt\label{balpha-tildabnalpha}
Let $\B_\a$ and $\tilde{\B}_{\a, n}$ be the unique solutions of \eqref{tikhexact2} and \eqref{finitetikhpert1opeq}, respectively. Let $C_0$ be as in Theorem \ref{dirichletestimate} and $\B$ be as in \eqref{exacteqn}. Let  $\tilde u$ be as in Assumption \ref{noise} and $\e_n>0$ be such that $\lVert{\T_{\tilde{u}}-\T_{\tilde{u}}\P_n}\rVert\leq \e_n$ and $\d>0$ be as given in \eqref{uestimate}. Then
$$\lVert{\B_\a-\tilde{\B}_{\a, n}}\rVert\leq  \dfrac{d\sqrt{C_0}\,\d}{\sqrt{\a}}\lVert{\B}\rVert +\dfrac{\e_n}{\sqrt{\a}}\lVert{\B}\rVert +\dfrac{\d}{2\sqrt{\a}}.$$
\et
\bpf
We have
\begin{eqnarray}
\B_\a-\tilde{\B}_{\a, n}&=&(\T_u^*\T_u+\a I)^{-1}\T_u^*(v_0-u)-(\P_n\T_{\tilde{u}}^* \T_{\tilde{u}}\P_n+\a I)^{-1}\P_n\T_{\tilde{u}}^*(v_0-\tilde{u})\nonumber\\
&=& [(\T_u^*\T_u+\a I)^{-1}\T_u^*-(\P_n\T_{\tilde{u}}^* \T_{\tilde{u}}\P_n+\a I)^{-1}\P_n\T_{\tilde{u}}^*](v_0-u)\label{s-1}\\
&&+(\P_n\T_{\tilde{u}}^* \T_{\tilde{u}}\P_n+\a I)^{-1}\P_n\T_{\tilde{u}}^*(\tilde{u}-u).\nonumber
\end{eqnarray}
Now, using \eqref{exacteqn}, we have
\beqarray
(\P_n\T_{\tilde{u}}^* \T_{\tilde{u}}\P_n+\a I)^{-1}\P_n\T_{\tilde{u}}^*(v_0-u)&=&(\P_n\T_{\tilde{u}}^* \T_{\tilde{u}}\P_n+\a I)^{-1}\P_n\T_{\tilde{u}}^*\T_u\B\\
(\T_u^*\T_u+\a I)^{-1}\T_u^*(v_0-u)&=&(\T_u^*\T_u+\a I)^{-1}\T_u^*\T_u\B.
\eeqarray
Also,
\begin{eqnarray}
(\P_n\T_{\tilde{u}}^* \T_{\tilde{u}}\P_n+\a I)^{-1}\P_n\T_{\tilde{u}}^*\T_u\B-(\T_u^*\T_u+\a I)^{-1}\T_u^*\T_u\B
&=&(\P_n\T_{\tilde{u}}^* \T_{\tilde{u}}\P_n+\a I)^{-1}\P_n\T_{\tilde{u}}^*
(\T_u\nonumber\\
&&-\T_{\tilde{u}}\P_n)\T_u^*\T_u(\T_u^*\T_u+\a I)^{-1}\B\label{s-2}\\
&&+\a(\P_n\T_{\tilde{u}}^* \T_{\tilde{u}}\P_n+\a I)^{-1}
(\P_n\T_{\tilde{u}}^*\nonumber\\
&&-\T_u^*)(\T_u\T_u^*+\a I)^{-1}\T_u\B.\nonumber
\end{eqnarray}
We now recall the following estimates from \cite{nairopeq},
\beqarray
\lVert{(\P_n\T_{\tilde{u}}^*\T_{\tilde{u}}\P_n+\a I)^{-1}\P_n\T_{\tilde{u}}^*}\rVert\leq \dfrac{1}{2\sqrt{\a}},&&\q \lVert{(\P_n\T_{\tilde{u}}^*\T_{\tilde{u}}\P_n+\a I)^{-1}}\rVert\leq \dfrac{1}{\a}\\
\lVert{\T_u^*\T_u(\T_u^*\T_u+\a I)^{-1}}\rVert\leq 1,&&\q \lVert{(\T_u\T_u^*+\a I)^{-1}\T_u}\rVert\leq \dfrac{1}{2\sqrt{\a}}\ .
\eeqarray
Using these estimates and using \eqref{s-2}, we obtain from \eqref{s-1},
\beqarray
\|\B_\a-\tilde{\B}_{\a, n}\| &\leq &\|[(\T_u^*\T_u+\a I)^{-1}\T_u^*-(\P_n\T_{\tilde{u}}^* \T_{\tilde{u}}\P_n+\a I)^{-1}\P_n\T_{\tilde{u}}^*](v_0-u)\|\\
&&+\lVert{(\P_n\T_{\tilde{u}}^* \T_{\tilde{u}}\P_n+\a I)^{-1}\P_n\T_{\tilde{u}}^*(\tilde{u}-u)}\rVert\\
 &\leq &\|{(\T_u^*\T_u+\a I)^{-1}\T_u^*\T_u\B
-(\P_n\T_{\tilde{u}}^* \T_{\tilde{u}}\P_n+\a I)^{-1}\P_n\T_{\tilde{u}}^*\T_u}\B\| \\
&&+\dfrac{\lVert{\tilde{u}-u}\rVert}{2\sqrt{\a}}\\
& \leq & \dfrac{\lVert{\T_u-\T_{\tilde{u}}\P_n}\rVert}{2\sqrt{\a}}\lVert{\B}\rVert +\dfrac{\lVert{\P_n\T_{\tilde{u}}^*-\T_u^*}\rVert}{2\sqrt{\a}}\lVert{\B}\rVert+\dfrac{\lVert{\tilde{u}-u}\rVert}{2\sqrt{\a}}.
\eeqarray
Now, using the fact $\|\P_n\T^*_{\tilde{u}}-\T^*_{u}\|=\|\T_u-\T_{\tilde{u}}\P_n\|$, we have
\beqarray
\|\B_\a-\tilde{\B}_{\a, n}\|&=&\dfrac{\lVert{\T_u-\T_{\tilde{u}}\P_n}\rVert}{\sqrt{\a}}\lVert{\B}\rVert+\dfrac{\lVert{\tilde{u}-u}\rVert}{2\sqrt{\a}}\\
&\leq &\dfrac{\lVert{\T_u-\T_{\tilde{u}}}\rVert}{\sqrt{\a}}\lVert{\B}\rVert+\dfrac{\lVert{\T_{\tilde{u}}-\T_{\tilde{u}}\P_n}\rVert}{\sqrt{\a}}\lVert{\B}\rVert+\dfrac{\lVert{\tilde{u}-u}\rVert}{2\sqrt{\a}}.
\eeqarray
\noindent
Therefore, using \eqref{uestimate} and \eqref{estimateoft}, we have
\beqarray
\lVert{\B_\a-\tilde{\B}_{\a, n}}\rVert &\leq & \dfrac{\lVert{\T_u-\T_{\tilde{u}}}\rVert}{\sqrt{\a}}\lVert{\B}\rVert+\dfrac{\lVert{\T_{\tilde{u}}-\T_{\tilde{u}}\P_n}\rVert}{\sqrt{\a}}\lVert{\B}\rVert + \dfrac{\lVert{\tilde{u}-u}\rVert}{2\sqrt{\a}}\\
&\leq & \dfrac{d\sqrt{C_0}\,\d}{\sqrt{\a}}\lVert{\B}\rVert +\dfrac{\e_n}{\sqrt{\a}}\lVert{\B}\rVert +\dfrac{\d}{2\sqrt{\a}}.
\eeqarray
\epf

Let  $$ \tilde{\A}_{\a, n}:=\tilde{\B}_{\a, n}+\A_0.$$ Then, our next theorem will give an estimate of $\|\A-\tilde{\A}_{n,\a}\|,$ the proof follows from Theorem \ref{balpha-tildabnalpha} and using the fact
$$\|\A-\tilde{\A}_{\a, n}\|= \|\B-\tilde{\B}_{\a, n}\|\leq\|\B-\B_\a\|+\|\B_\a-\tilde{\B}_{\a, n}\|\ .$$

\bt\label{final}
Let $\B_\a$ and $\tilde{\B}_{\a,n}$ be the unique solutions of \eqref{tikhexact2} and \eqref{finitetikhpert1opeq}, respectively. Let $C_0$ be as in Theorem \ref{dirichletestimate} and $\B$ be as in \eqref{exacteqn}. Let $\e_n>0$ be such that $\lVert{\T_{\tilde{u}}-\T_{\tilde{u}}\P_n}\rVert\leq \e_n$ and $\d>0$ be as given in \eqref{uestimate}. Then
$$\|\A-\tilde{\A}_{\a,n}\| \leq \lVert{\B-\B_\a}\rVert +\dfrac{d\sqrt{C_0}\d}{\sqrt{\a}}\lVert{\B}\rVert +\dfrac{\e_n}{\sqrt{\a}}\lVert{\B}\rVert +\dfrac{\d}{2\sqrt{\a}}.$$
\et
\brem\label{epsilonto0}
By Corollary \ref{cpt-appr}, we know that   $\lVert{\T_{\tilde{u}}-\T_{\tilde{u}}\P_n}\rVert\to 0$ as $n\to \infty.$ Hence $\e_n>0$ can be chosen in such a way that $\lVert{\T_{\tilde{u}}-\T_{\tilde{u}}\P_n}\rVert\leq \e_n$ and $\e_n\to 0$ as $n\to \infty$. Thus, Theorem \ref{final} shows that,  by an appropriate choice of  $\a$ and $n$, depending on  $\d$,   $\A_{\a, n}$ is an approximation of $\A$. 
\erem

\section{Adaptive choice of Parameters}\label{adaptive}
Recall that, our aim is to show $\tilde{\A}_{\a,n}$ is a stable approximation to  $\A$ for some chosen parameters $n$ and $\a$ compatible with the noise level $\d.$ In other words, we need to choose the parameters $n$ and $\a$, depending on $\d$, in such a way that $\tilde{\A}_{\a,n}$ converges to $\A$ as the noise level $\d\to 0.$ In this regard, we observe that Theorem \ref{final} will show us a direction for choosing the parameters suitably. Parameter choice strategy in regularization theory has a vast literature, for more details one may look into \cite{neuber}, \cite{nairopeq}, \cite{schock} and the references therein. In this paper we will use the procedure adopted in \cite{schock}.

Now, we would like to recall that $\B=\A-\A_0$ and $\B_\a$ is the solution of \eqref{tikhexact2}, that is, a Tikhonov regularized solution. Therefore, from the theory of Tikhonov regularization, it is known that (cf. \cite{nairopeq}) $$\|\B-\B_\a\|\to 0\q \text{as}\ \a\to 0.$$
But, for obtaining an estimate for $\|\B-\B_\a\|$, it is necessary to assume some source condition on $\B.$ So, we make the following general assumption for source condition.

\textbf{Source condition:} Let
\beq\label{source1}
\B=\varphi(\T_u^*\T_u)\C,\q \|\C\|\leq \r
\eeq
for some $\r>0$ and for some monotonically increasing function $\varphi$ defined on $(0, \g],$ where $\g\geq \|\T_u\|^2,$ such that $\ds\lim_{\lambda\to 0}\varphi(\lambda)=0$ and $$\ds \sup_{0\leq \lambda\leq \g}\dfrac{\a\varphi(\lambda)}{\lambda+\a}\leq \varphi(\a)\q \text{for all}\ \a>0.$$
A typical case of such a situation is when $\B$ is in the range of $\f(\T_u^*\T_u)$, where
$\f(\l):=\l^\nu$  for some
$\nu\in (0,1]$ or
$\f(\l):=[\log(1/\l)]^{-p}$ for some $p>0$ (see, for examaple, \cite{neuber, nairopeq}).

Using \eqref{tikhpertbd2}, \eqref{tikhexact2} and the above assumption, we have
$$\B-\B_\a=\a(\T_u^*\T_u+\a I)^{-1}\B=\a(\T_u^*\T_u+\a I)^{-1}\varphi(\T_u^*\T_u)\C$$
so that
\beqarray
\|\B-\B_\a\| &=& \|\a(\T_u^*\T_u+\a I)^{-1}\varphi(\T_u^*\T_u)\C\|\\
&\leq & \sup_{0\leq \lambda\leq \g}\dfrac{\a\varphi(\lambda)}{\lambda+\a} \|\C\|\leq \r \varphi(\a).
\eeqarray
Thus, with the above source condition, and by Theorem \ref{final}, we obtain
\beq\label{finalestimate}
\|\B-\tilde{\B}_{\a, n}\|\leq \r\varphi(\a)+C_{\B} \dfrac{ \d(d+1)+\e_n}{\sqrt{\a}},
\eeq
where $C_{\B}=\max\Big\{1/2, \sqrt{C_0}\,\|\B\|\Big\}.$

As mentioned in Remark \ref{epsilonto0}, $\e_n\to 0$ as $n\to \infty.$ Let $\ds n_\d\in \N$ be such that $\e_n<\d$ for all $n\geq n_\d.$ Then by \eqref{finalestimate}, we have
\beq\label{deltaestimate}
\|\B-\tilde{\B}_{\a, n_\d}\|\leq \r\varphi(\a)+C_{\B} \dfrac{ \d(d+2)}{\sqrt{\a}}.
\eeq
Let $C>1/2$ be a constant such that $\r \varphi(\g)\leq C.$ Then, using \eqref{source1}, we obtain $\|\B\|\leq C$, and hence $$C_{\B}=\max\Big\{1/2,\sqrt{C_0}\, \|\B\|\Big\}\leq C.$$
Thus, by \eqref{deltaestimate} we have
\beq\label{cdeltaestimate}
\|\B-\tilde{\B}_{\a, n_\d}\|\leq \r\varphi(\a)+C\, \dfrac{ \d(d+2)}{\sqrt{\a}}.
\eeq

We now apply the adaptive technique for choosing the parameter $\a$ \textit{a posteriori}, so that the rate of convergence is order optimal. Here we would like to point out that, the technique is elaborately explained in general setting in \cite{schock}. Also, whatever result we will obtain, can be derived easily from the results already obtained in \cite{santosh}, \cite{schock}. But, in order to keep this paper self contained we are giving the details, presented in our own way.

Let $\d>0$ be such that $\d^2(d+2)^2\leq \g.$ Let $\a_0=\d^2(d+2)^2.$ Since $\varphi$ is an increasing function on $(0,\g]$ and $\a_0\leq \g,$ we have
$$\r\varphi(\a_0)\leq \r\varphi(\g)\leq C.$$
Let $\mu >1$ be any fixed real number and $N\in \N$ be fixed. Let
\beq\label{alphai}
\a_i=\mu^{2i}\a_0,\q i=1,2,...,N.
\eeq
Then, clearly we have $$0<\a_0<\a_1<...<\a_N$$and
$$\sqrt{\a_i}\leq \mu \sqrt{\a_{i-1}}\q \text{for all}\ i=1,2,...,N.$$
Let \beq\label{l}
l=\max\Big\{i\in \{0,1,...,N-1\}:\r\mu^i\varphi(\a_i)\leq C\Big\}.
\eeq
We now establish the following lemma.
\bl\label{l-j}
Let $n_\d$ and $l$ be as in \eqref{deltaestimate} and \eqref{l}, respectively. Then, for any $j\in \{0,1,...,l\},$ we have
$$\|\tilde{\B}_{\a_l, n_\d}-\tilde{\B}_{\a_j, n_\d}\|\leq \dfrac{4C}{\mu^j}.$$
\el
\bpf
First we note that for any $j\leq l,$ we have $$\r\mu^j\varphi(\a_j)\leq \r\mu^l\varphi(\a_l).$$
Then, using \eqref{cdeltaestimate}, \eqref{alphai} and \eqref{l}, we have
\beqarray
\|\tilde{\B}_{\a_l, n_\d}-\tilde{\B}_{\a_j, n_\d}\| &\leq &\|\tilde{\B}_{n_\d,\a_l}-\B\| + \|\tilde{\B}_{n_\d,\a_j}-\B\|\\
&\leq & \r\varphi(\a_l)+\dfrac{C}{\mu^l}+\r\varphi(\a_j)+\dfrac{C}{\mu^j}\\
&\leq & \dfrac{4C}{\mu^j}.
\eeqarray
\epf
Let \beq\label{k}
k=\max\Big \{i\in \{0,1,...,N\}:\|\tilde{\B}_{n_\d,\a_i}-\tilde{\B}_{n_\d,\a_j}\|\leq \dfrac{4C}{\mu^j}, j=0,1,...,i\Big\}.
\eeq
Then, Lemma \ref{l-j} ensures that $l\leq k.$ Now, using \eqref{cdeltaestimate}, \eqref{l}, \eqref{k} we have
\beqarray
\|\B-\tilde{\B}_{\a_k, n_\d }\| &\leq & \|\B- \tilde{\B}_{\a_l, n_\d}\|+\|\tilde{\B}_{\a_l, n_\d}-\tilde{\B}_{\a_k, n_\d}\|\\
&\leq & \r\varphi(\a_l)+\dfrac{C}{\mu^l}+\dfrac{4C}{\mu^l}\\
&\leq & \dfrac{2C}{\mu^l} +\dfrac{4C}{\mu^l}\\
&=&\dfrac{6C}{\mu^l}.
\eeqarray
Let  $\a_\d$ be such that $\r\varphi(\a_\d)=C\dfrac{\d(d+2)}{\sqrt{\a_\d}}.$ Then, it is clear that $\r\varphi(\a)+C\dfrac{\d(d+2)}{\sqrt{\a}}$ attains it minimum at $\a_\d.$ Now, using the definition of $l,$ we have
\beqarray
\varphi(\a_\d)\sqrt{\a_\d}=C\dfrac{\d(d+2)}{\r}<\d(d+2)\mu^{l+1}\varphi(\a_{l+1})&=&\sqrt{\d^2(d+2)^2\mu^{2(l+1)}}\ \varphi(\a_{l+1})\\
&=&\sqrt{\a_{l+1}}\ \varphi(\a_{l+1}).\eeqarray
Now, using the fact that $\varphi$ is an increasing function, we obtain$$\a_\d<\a_{l+1}.$$
Thus, we have $$\sqrt{\a_\d}<\sqrt{\a_{l+1}}=\sqrt{\a_0}\mu^{l+1}=\d(d+2)\mu^{l+1}.$$
Therefore,
\beq\label{optimal}
\|\B-\tilde{\B}_{\a_k, n_\d}\|\leq \dfrac{6C}{\mu^l}\leq \dfrac{6C\mu\d(d+2)}{\sqrt{\a_\d}}=6\mu \r\varphi(\a_\d).\eeq
Let \beq\label{psi}
\Psi(\lambda)=\dfrac{\r\lambda\sqrt{\varphi^{-1}(\lambda)}}{C(d+2)},\q 0<\lambda\leq \g.
\eeq
Then, we have$$\d=\dfrac{\r \varphi(\a_\d)\sqrt{\a_\d}}{C(d+2)}=\Psi(\varphi(\a_\d)),$$
so that $\varphi(\a_\d)=\Psi^{-1}(\d).$
Thus, using \eqref{optimal}, we have
\beq\label{opt}
\|\B-\tilde{\B}_{\a_k, n_\d}\|\leq 6\mu \r\Psi^{-1}(\d)\ .
\eeq
Since $\B=\A-\A_0,$ therefore, using \eqref{opt}, we have obtained the following theorem .
%about optimal rate of convergence for \textit{a posteriori} parameter choice.
\bt\label{optimalconvergence}
Let $\r,\mu,k$ be as defined in \eqref{source1}, \eqref{alphai} and \eqref{k}, respectively and $\A_0$ be as in \eqref{v0}. Let $n_\d$ and $\a_k$ be as defined in \eqref{deltaestimate} and \eqref{alphai}, respectively, and let $\Psi$ be as defined in \eqref{psi}. Let $\tilde{\A}_{\a_k, n_\d}:=\A_0+\tilde{\B}_{\a_k, n_\d}.$ Then, we have
$$\|\A-\tilde{\A}_{\a_k, n_\d}\|\leq 6\mu \r\Psi^{-1}(\d)\ .$$
\et
 The above theorem ensures that the \textit{a posteriori} chosen parameter $\a_k$ satisfies the optimal rate of convergence with respect to $\d$, without the knowledge of any \textit{a priori} source function. So, the procedure of adaptive choosing of the regularizing parameter is effective in our analysis of obtaining the \textit{a posteriori} regularizing parameter.

\section{Smoothing of noisy data}\label{smoothingnoisy}
Recall that in Remark \ref{smoothingremark}, we have mentioned the need for smoothed version of the noisy data $\tilde{u}\in L^2(0,\tau;L^2(\O)).$ This section is devoted to that purpose. In \cite{samprita}, the authors have used the Clement operator (see \cite{clement}) for obtaining a smoothed version of a noisy data in the context of a parameter identification problem for an elliptic PDE. Here we will do similar kind of smoothing, but we will be doing for the parabolic case. 

Let $u$ be as considered in the Assumption \ref{existence}. Then $u\in L^2(0,\tau;L^2(\O))$ along with $\nabla u\in L^2(0,\tau;L^\infty(\O,\R^d)).$ In order to do our smoothing analysis, we assume that $u\in L^2(0,\tau; H^4(\O)).$ Let $\tilde{u}\in L^2(0,\tau;L^2(\O))$ be the observed data corresponding to $u.$ We want to obtain an element, say $z\in L^2(0,\tau;L^2(\O))$ such that $\nabla z\in L^2(0,\tau;L^\infty(\O,\R^d))$, and we will call this $z$ as a smoothed version of $\tilde{u}.$ For this we assume that $\O$ is a polygonal domain in $\R^2$.

Let $\mathcal{L}$ be an element of a family of quasi uniform triangulation of $\O.$ Then there exist a constant $\g_0>0$ such that
 \beq\label{quasitriang}
  \min_{S\in \mathcal{L}}\dfrac{\text{diam}(S)}{h}\geq \g_0>0,
  \eeq where $\text{diam}(S)$ is the diameter of triangle $S\in \mathcal{L}$ and $h=\max_{S\in \mathcal{L}}\text{diam}(S),$ is the mesh size.
Let $\Pi$ be the Clement operator (see \cite{clement}) that maps $L^2(\O)$ to the space of all polynomials of degree less than or equal to $3$. We now state an important result required for our analysis, for its proof refer to \cite{clement}.

\bt{\rm(\cite{clement}, Theorem 1)}\label{1}
Let $v\in L^2(\O)$. Then $\Pi v\in W^{1,\infty}(\O)$ and  there exist constants $C_1, C_2>0$ such that $$\|v-\Pi v\|_{L^2}\leq C_1\|v\|_{L^2}\q \text{for all}\ v\in L^2(\O)$$and$$\|v-\Pi v\|_{H^3}\leq C_2 h\|v\|_{H^4}\q \text{for all}\ v\in H^4(\O).$$
\et

We know that if $\O\subset \R^d$ with sufficiently smooth boundary, then for every $k$ with $d<2k$, we have the continuous embedding $$H^k(\O)\hookrightarrow L^\infty(\O).$$
In our case, since $u\in L^2(0,\tau;H^4(\O))$ and $\O\subset \R^2$, therefore we have the continuous embedding $$H^4(\O)\hookrightarrow H^3(\O)\hookrightarrow W^{1,\infty}(\O),$$ and hence there exist a constant $C_3>0$ such that
\beq\label{embd}
\|u(\cdot, t)-\Pi u(\cdot, t)\|_{W^{1,\infty}}\leq C_3\| u(\cdot, t)-\Pi u(\cdot, t)\|_{H^3}\q \text{for a.a}\ t\in [0,\tau].\eeq
Let $C_4=C_2C_3$. Then, using Theorem \ref{1} and \eqref{embd}, we have
\beq\label{c4}
\|u(\cdot, t)-\Pi u(\cdot, t)\|_{W^{1,\infty}}\leq C_4h\|u(\cdot, t)\|_{H^4}\q \text{for a.a}\ t\in [0,\tau].
\eeq
We now state a result, which is a reformulation of Theorem 4.5.11 in \cite{brenner}.
\bl
For any $S\in \mathcal{L}$ and for all $v\in L^2(S)$, we have
$$\|\Pi v\|_{W^{1,\infty}(S)}\leq \dfrac{1}{(\text{diam}(S))^2}\|\Pi v\|_{{L^2(S)}}.$$
\el
\noindent Let $v\in L^2(\O).$ Then by the above Lemma, we have
\beqarray
\|\Pi v\|_{W^{1,\infty}(S)}& \leq & \dfrac{1}{(\text{diam}(S))^2}\|\Pi v\|_{L^2(S)}\leq \dfrac{1}{h^2\g^2_0}\|\Pi v\|_{L^2(S)}\\
&\leq & \dfrac{1}{h^2\g^2_0}\|\Pi v\|_{L^2(\O)}\leq \dfrac{1}{h^2\g^2_0}\big(\|(1-\Pi)v\|_{L^2(\O)}+\|v\|_{L^2(\O)}\big)\\
&\leq & \dfrac{C_1+1}{h^2\g^2_0}\|v\|_{L^2(\O)},
\eeqarray
where $C_1$ is as in Theorem \ref{1} and $\g_0$ is as in \eqref{quasitriang}. Thus, using the fact that $$\|\Pi v\|_{W^{1,\infty}(\O)}\leq \max_{S\in \mathcal{L}}\|\Pi v\|_{W^{1,\infty}(S)}$$ we have proved the following theorem.
\bt\label{c5}
 Let $C_5:=\dfrac{C_1+1}{\g^2_0}$, where $C_1$ is as in Theorem \ref{1} and $\g_0$ is as in \eqref{quasitriang}. Then for all $v\in L^2(\O),$
  $$\|\Pi v\|_{W^{1,\infty}(\O)}\leq \dfrac{C_5}{{h^2}}\|v\|_{L^2(\O)}.$$
\et
\section{Modified error estimates}
Let $\Pi$ be as considered in Section \ref{smoothingnoisy} and $u$ be as considered in Assumption \ref{existence}. Let the noisy data $\tilde{u}\in L^2(0,\tau;L^2(\O))$ be such that
 \beq\label{delta}
 \|u-\tilde{u}\|_{L^2(0,\tau;L^2(\O))}\leq \d.
 \eeq
Using the results of Section \ref{smoothingnoisy}, we now obtain error estimates with the noisy data $\tilde{u}.$
\bt\label{dh}
Let $\d>0$ be as considered in \eqref{delta}. Let $C_4$ and $C_5$ be as in \eqref{c4} and Theorem \ref{c5}, respectively. Let $C=max\{C_4,C_5\}.$ Then
$$\Big(\int_0^\tau\|\nabla u(\cdot, t)-\nabla\Pi \tilde{u}(\cdot, t)\|^2_{L^\infty}dt\Big)^{1/2}\leq C\Big(h\|u\|_{L^2(0,\tau;H^4(\O))}+\dfrac{\d}{{ h^2}}\Big).$$
\et
\bpf
Using \eqref{c4} and Theorem \ref{c5}, we have
\beqarray
\|u(\cdot, t)-\Pi \tilde{u}(\cdot, t)\|_{W^{1,\infty}}&\leq & \|u(\cdot, t)-\Pi u(\cdot, t)\|_{W^{1,\infty}}+\|\Pi u(\cdot, t)-\Pi \tilde{u}(\cdot, t)\|_{W^{1,\infty}}\\
&\leq & C_4h\|u(\cdot, t)\|_{H^4}+\dfrac{C_5}{h^2}\|u(\cdot, t)-\tilde{u}(\cdot, t)\|_{L^2}
\eeqarray
for a.a $t\in [0,\tau].$ Thus, we have
\beqarray
\int_0^\tau\|\nabla u(\cdot, t)-\nabla\Pi \tilde{u}(\cdot, t)\|^2_{L^\infty}dt &\leq & \int_0^\tau\| u(\cdot, t)-\Pi \tilde{u}(\cdot, t)\|^2_{W^{1,\infty}}dt\\
&=&C_4^2h^2\|u\|^2_{L^2(0,\tau;H^4(\O))}+\dfrac{C_5^2}{h^4}\int_0^\tau\|u(\cdot, t)-\tilde{u}(\cdot, t)\|^2_{L^2}dt\\
&&+\dfrac{2C_4C_5}{h}\int_0^\tau\|u(\cdot, t)\|_{H^4}\|u(\cdot, t)-\tilde{u}(\cdot, t)\|_{L^2}dt\\
&\leq &C_4^2h^2\|u\|^2_{L^2(0,\tau;H^4(\O))}+\dfrac{C_5^2}{h^4}\|u-\tilde{u}\|^2_{L^2(0,\tau;L^2(\O))}\\
&&+\dfrac{2C_4C_5}{h}\|u\|_{L^2(0,\tau;H^4(\O))}\|u-\tilde{u}\|_{L^2(0,\tau;L^2(\O))}\\
&=&\Big(C_4h\|u\|_{L^2(0,\tau;H^4(\O))}+\dfrac{C_5}{h^2}\|u-\tilde{u}\|_{L^2(0,\tau;L^2(\O))}\Big)^2
\eeqarray
and hence
$$
\Big(\int_0^\tau\|\nabla u(\cdot, t)-\nabla\Pi \tilde{u}(\cdot, t)\|^2_{L^\infty}dt\Big)^{1/2}\leq C\Big(h\|u\|_{L^2(0,\tau;H^4(\O))}+\dfrac{\d}{h^2}\Big).
$$
\epf
\noindent {We define $\widehat{\Pi}:L^2(0,\tau;L^2(\O))\to L^2(0,\tau;L^2(\O))$ by $$(\widehat{\Pi}\psi)(t)(\cdot):=\Pi( \psi(\cdot,t))\q \text{for }\ t\in [0,\tau].$$
With this definition of $\widehat{\Pi}$ we have the following theorem.}
\bt\label{u-pi}
{Let $\d>0$ be as considered in \eqref{delta}. Let $C_4$ and $C_5$ be as in \eqref{c4} and Theorem \ref{c5}, respectively. Let $C=max\{C_4,C_5\}$ and $|\O|$ denotes the Lebesgue measure of the set $\O$ and $C'=C\sqrt{|\O|}.$ Then $$\|u-\widehat{\Pi} \tilde{u}\|_{L^2(0,\tau;L^2(\O))}\leq C'\Big(h\|u\|_{L^2(0,\tau;H^4(\O))}+\dfrac{\d}{h^2}\Big).$$}
\et
\bpf
First we observe that $$\|u-\widehat{\Pi} \tilde{u}\|^2_{L^2(0,\tau;L^2(\O))}=\int_0^\tau\|u(\cdot,t)-\Pi \tilde{u}(\cdot,t)\|^2_{L^2(\O)}dt\leq |\O|\int_0^\tau\|u(\cdot, t)-\Pi \tilde{u}(\cdot, t)\|^2_{W^{1,\infty}(\O)}dt.$$
Now, following the proof of Theorem \ref{dh}, we obtain
$$\int_0^\tau\|u(\cdot, t)-\Pi \tilde{u}(\cdot, t)\|^2_{W^{1,\infty}(\O)}dt\leq C^2\Big(h\|u\|_{L^2(0,\tau;H^4(\O))}+\dfrac{\d}{h^2}\Big)^2.$$
Thus, we have$$\|u-\widehat{\Pi} \tilde{u}\|_{L^2(0,\tau;L^2(\O))}\leq C'\Big(h\|u\|_{L^2(0,\tau;H^4(\O))}+\dfrac{\d}{h^2}\Big).$$
\epf

By Theorem \ref{dh} and Theorem \ref{u-pi}, we have 
$$\|u-\widehat{\Pi}\tilde{u}\|_{L^2(0,\tau;L^2(\O))}+\Big(\int_0^\tau\|\nabla u(\cdot,t)-\nabla\Pi \tilde{u}(\cdot,t)\|^2_{L^\infty}dt\Big)^{1/2}\leq  2\tilde{C}\d_h,$$
where 
\beq\label{constant} \tilde{C}=\max\Big\{C\|u\|_{L^2(0,\tau;L^2(\O))},\ C'\|u\|_{L^2(0,\tau;L^2(\O))}\Big\}\eeq and
\beq\label{deltah}
\d_h:= \max\Big\{h,\dfrac{\d}{h^2}\Big\}.
\eeq
Now, corresponding to the Assumption \ref{noise} on $\tilde u$, we have the inequality 
$$\|u-\widehat{\Pi}\tilde{u}\|_{L^2(0,\tau;L^2(\O))}+\Big(\int_0^\tau\|\nabla u(\cdot,t)-\nabla\Pi \tilde{u}(\cdot,t)\|^2_{L^\infty}dt\Big)^{1/2}\leq  2\tilde{C}\d_h,$$
for $\widehat{\Pi} \tilde{u}$, where $\tilde{C}$ is as in (\ref{constant}).
Then carrying out the analysis as done in Section \ref{reg} and Section \ref{finit}, with $\tilde{u}$ and $\d$ replaced by $\widehat{\Pi} \tilde{u}$ and $\d_h$, respectively, by Theorem \ref{final}, we obtain the following theorem.
\bt
Let $\B_\a$ be the unique solution of \eqref{tikhexact2} and $\widehat{\B}_{\a, n}$ be the unique solution of \eqref{finitetikhpert1opeq}, with $\tilde{u}$ replaced by $\widehat{\Pi} \tilde{u}$. Let $C_0$ be as in Theorem \ref{dirichletestimate}. Let $\B$ be as in \eqref{exacteqn} and $\widehat{\A}_{\a, n}:=\widehat{\B}_{\a, n}+\A_0.$ Let $\tilde{\e}_n>0$ be such that $\lVert{\T_{\widehat{\Pi} \tilde{u}}-\T_{\widehat{\Pi}\tilde{u}}\P_n}\rVert\leq \tilde{\e}_n$ and $\d_h>0$ be as given in \eqref{deltah}. Then
\beq\label{deltahestimate}
\|\A-\widehat{\A}_{\a, n}\|\leq \lVert{\B-\B_\a}\rVert +C''\,\dfrac{\d_h+\tilde{\e}_n}{\sqrt{\a}}
\eeq
where $\|\B-\B_\a\|\to 0$ as $\a\to 0$ and $C''$ is a positive constant depending only on the constants $\|\B\|, \|u\|_{L^2(0,\tau;L^2(\O))}, d, \tilde{C}$ and $\sqrt{C_0}.$
\et

\brem Note that,  if we fix the mesh size $h$ first, and then choose the error level $\d$ in such a way that $\d\leq h^3,$ then  $\d_h=h$ and hence the estimate in \eqref{deltahestimate} becomes$$\|\A-\widehat{\A}_{\a, n}\|\leq \lVert{\B-\B_\a}\rVert +O\Big(\dfrac{h+\tilde{\e}_n}{\sqrt{\a}}\Big).$$
\erem

\section{Conclusion}
We have considered an inverse problem of identifying a coefficient $\A\in (H^1(\O))^{d\times d}$ of a parabolic PDE with Dirichlet boundary condition. Under specific assumption, a uniqueness result for the solution of inverse problem is obtained.  By making use of a weak formulation, we have reduced our inverse problem into solving an ill-posed operator equation, where the operator involved is linear, and we have explicitly obtained a representation for the adjoint of the corresponding linear operator. In order to obtain stable approximations for the sought coefficient $\A,$ we have used the theory of Tikhonov regularization. Also, we have given a finite dimensional realization of the method for practical implementation. For the parameter choice, we have used the adaptive technique to obtain the regularizing parameter effectively. Finally, we have demonstrated a procedure to obtain a smoothed version of a noisy data by making use of Clement operator. But, we would like to mention that for smoothing, we have assumed a higher regularity of the data, namely $u\in L^2(0,\tau;H^4(\O)).$

\end{document}